\documentclass[10pt]{article}
\usepackage{latexsym}
\usepackage{verbatim}
\usepackage{enumerate}
\usepackage[dvips,draft,final]{graphics}

\newtheorem{defi}{Definition}
\newtheorem{lem}[defi]{Lemma}
\newtheorem{prop}[defi]{Proposition}

\newtheorem{theo}[defi]{Theorem}

\newtheorem{rmq}[defi]{Remark}

\newenvironment{dem}{{\bf Proof:}}{ { \hspace {\stretch{1} }$\Box$}}

\newcommand{\Tr}{\mathrm{Tr}}

\newcommand{\N}{{\mathbb N}}

\newcommand{\R}{{\mathbb R}}
\newcommand{\C}{{\mathbb C}}
\newcommand{\1}{\underline{1}}

\newcommand{\Hh}{{\mathcal H}}

\newcommand{\eps}{\varepsilon}
\def\C{{\bf C}}

\def\N{{\bf N}}
\def\R{{\bf R}}

\def\CC{{\mathcal C}}

\def\eps{\varepsilon}

\def\ni{\noindent}
\def\op{\op}

\def\ol{\overline}
\def\d{{\rm d}}

\def\op_#1{\mathrel{\mathop{{\rm op}_{#1}}}}

\def\build#1_#2^#3{\mathrel{
\mathop{\kern 0pt#1}\limits_{#2}^{#3}}}

\def\td_#1,#2{\mathrel{
\mathop{\build\longrightarrow_{#1\rightarrow #2}^{}}}}

\def\limsup_#1,#2{\mathrel{
\mathop{\build{\rm limsup}_{#1\rightarrow#2}^{}}}}

\def\liminf_#1,#2{\mathrel{
\mathop{\build{\rm liminf}_{#1\rightarrow#2}^{}}}}

\def\aref#1{(\ref{#1})}
\def\eps{\varepsilon}

\def\1{{\bf 1}}
\def\0{{\bf 0}}

\def\ni{\noindent}
\def\ol{\overline}
\def\d{{\rm d}}

\begin{document}

\title{Semi-classical analysis of a conjoint crossing of three symmetric modes.}

\author{Clotilde FERMANIAN KAMMERER\\ {\it \small Universit\'e Paris 12 -- UMR 8050} \\
Vidian ROUSSE \\
{\it \small Universit\'e Paris 12  -- UMR 8050}}

\date{ }

\maketitle

{\small{\bf Abstract}: {\it  In this article we focus on a semiclassical Schr\"odinger equation with matrix-valued potential presenting
 a symmetric conjoint crossing of three eigenvalues. The potential we consider  is well-known in the chemical literature as a  pseudo
 Jahn-Teller potential. We analyze the energy transfers which occur between the three modes in terms of Wigner measures.}}

\section{Introduction}

 We consider the following Schr\"odinger equation
\begin{equation}\label{eq:schro}
\left\{\begin{array}{l} i\eps\partial_t
\psi^\eps(q,t)=\left(-{\eps^2\over
2}\Delta_q+V(q)\right)\psi^\eps(q,t),\;\;(q,t)\in\R^2\times\R\\*[1ex]
\psi^\eps(q,0)=\psi^\eps_0(q)\in L^2(\R^2,\C^3),
\end{array}
\right.
\end{equation} where $\eps>0$ is the semi-classical  parameter and where the matrix-valued potential $V(q)=V_{PJT}(q)$ is
\begin{equation}\label{def:potential}
\forall q\in\R^2,\;\;V_{PJT}(q)=\left(\begin{array}{ccc} q_1 & 0 & q_2/\sqrt{2} \\ 0 & -
q_1 & q_2/\sqrt{2} \\ q_2/\sqrt{2} & q_2/\sqrt{2} & 0
\end{array}\right).
\end{equation}

\noindent  Therefore, the potential $V$ is subquadratic and the operator $-{\eps^2\over 2} \Delta_q+ V(q)$ is essentially self--adjoint on the Schwartz space.  Such systems naturally arise in quantum chemistry when
studying the
dynamics of large molecules in the frame of Born-Oppenheimer approximation. The parameter $\eps$ is  small
because it encounters for the ratio of the mass of the electron and of
the averaged mass of the nuclei (see  \cite{bo} or \cite{Ha94} for
the derivation of such equations). The solution $\psi^\eps(t) $
 does not have any direct
physical interpretation but only quadratic functions of it: for
example, for scalar equations, the position density $|\psi^\eps(q,t)|^2$  gives the
probability of finding the nuclei in the configuration $q\in\R^2$
at time $t$. Therefore, one is interested in the asymptotic
description for the time evolution of the matrix-valued Wigner function
of~$\psi^\eps(q,t)$,
$$W^\eps\!\left(\psi^\eps(t)\right)\!(q,p)=
(2\pi)^{-d}\int_{\R^2} \psi^\eps\!\left(q-{\textstyle{\eps\over
2}}v,t\right)\otimes \ol \psi^\eps\!\left(q+{\textstyle{\eps\over
2}}v,t\right)\,{\rm e}^{i\,v\cdot p}\,\d v$$ which plays the role
of a generalized probability density on phase space. Indeed, the
action of the Wigner function against compactly supported smooth
test functions $a\in\CC^\infty_c(\R^{4},\C^{3\times
3})$ is simply expressed in terms of the semi-classical
pseudodifferential operator with symbol $a$, which is defined for $\psi\in L^2(\R^2,\C^3)$ by
$$
\op_\eps(a)\psi(q)=(2\pi\eps)^{-d} \int_{\R^{4}}
a\left(\frac{q+v}{2},p\right){\rm
e}^{{i\over\eps}p\cdot(q-v)}\psi(v)\,\d v\,\d p.
$$
We then  have $\displaystyle{ {\rm tr}
\int_{\R^{4}}\,W^\eps(\psi)(q,p)a(q,p)\,\d q\,\d
p=\left(\op_\eps(a)\psi\,,\,\psi\right)_{L^2(\R^2,\C^3)}. }$ For
example, one recovers the position density by $\displaystyle{
|\psi(q)|^2={\rm tr} \int_{\R^2}  W^\eps(\psi)(q,p)\,\d p. }$ We
are concerned here with the description of the weak limits of the
Wigner transform, these distributions are  matrix-valued measure
called Wigner measure or semi-classical measures (see \cite{LP},
\cite{GL} or \cite{GMMP}).  Besides, these matrix-valued measures $\mu=(\mu_{i,j})$  are positive in the sense that their diagonal elements $\mu_{i,i}$ are positive Radon measures and their off-diagonal elements $\mu_{i,j}$ for $i\not=j$ are absolutely continuous with respect to $\mu_{i,i}$ and $\mu_{j,j}$.

$ $

 The potential $V_{PJT}(q)$ arises in models of quantum chemistry: it
is the  simplest pseudo Jahn-Teller Hamiltonian as introduced in
\cite{DYK} Chapter 10 (see also \cite{Viel} or \cite{Eisfeld} for
example).  Its main feature is that its eigenvalues  are
$\displaystyle{\sqrt{q_1^2+q_2^2}, \;\;-\sqrt{q_1^2+q_2^2}}$ and $0$. They  are
symmetric: the eigenvalue~$0$ is exactly the half sum of the two other eigenvalues. Besides, they are of
multiplicity $1$ as long as $(q_1,q_2)\not=(0,0)$ and  they simultaneously cross on the point $\{q_1=q_2=0\}$. We will
say that $V$ presents a  {\it conjoint symmetric crossing} in $q=(0,0)=0$. The appellation of {\it pseudo} Jahn-Teller
potential is a reference to what is called Jahn-Teller potentials
in quantum chemistry, namely the $2$ by $2$ potential of the form
${\displaystyle{\left(\begin{array} {cc}q_1 & q_2 \\ q_2 &
-q_1\end{array}\right)}}$ to which numerous mathematical works
have been devoted. Then, another simplest realization of  a
symmetric conjoint crossing is the matrix
$$ V_{JT}(q_1,q_2)= \left(\begin{array} {ccc}q_1 & q_2 & 0 \\ q_2 & -q_1 & 0 \\ 0 & 0 & 0 \end{array}\right).$$
Both matrices $V_{JT}$ and $V_{PJT}$ (defined in \aref{def:potential})
have the same symmetric eigenvalues and the behavior of solutions
to \aref{eq:schro} with the Jahn-Teller potential $V_{JT}$ is
well understood. Indeed, in that situation, the third mode evolves
at leading order independently from the two other ones which
interacts following the well-known process of conical
intersections (see \cite{La}, \cite{Ze}, \cite{Ha94}, \cite{HJ},
\cite{FG03} for example). However, the situation is different  for
the pseudo Jahn-Teller potential $V_{PJT}(q)$, we will see that the
three modes will interact altogether; note also that for this
potential, the three spectral projectors are singular on
$\{q_1=q_2=0\}$ while $V_{JT}$ has a smooth eigenprojector (for the mode~$0$).
While an important literature have been devoted to two by two crossings, there are very few results on conjoint crossing of three eigenvalues. The aim of this work is to analyze on a simple model the mechanism of a symmetric crossing of three eigenvalues.

$ $

We shall describe how the three modes
interact together for the pseudo Jahn-Teller potential $V_{PJT}(q)$ by
studying the branching of Wigner measures on the
crossing set $\{q_1=q_2=0\}$. We explain our result in the
following section. The main ingredient of the proof are a normal
form result given in Section~\ref{sec:nf} and a scattering theorem in Section~\ref{sec:scat}. The normal form reduces the analysis
of the crossing  to the study of
   a system  very close
to the system of three ordinary differential equations
\begin{equation}\label{eq:1}
i\eps\partial_s u + \left(\begin{array} {ccc} s & 0 & {z/\sqrt 2}\\  0 & -s & z/\sqrt 2\\ z/\sqrt2 & z/\sqrt 2 & 0 \end{array}\right)=0  .
\end{equation}
The behavior of solutions of this system  as $\eps$ goes to $0$ is analyzed in 
Section~\ref{sec:scat}. This system has been already studied by physicists (\cite{BE} and
\cite{CH}); with their arguments, one can find the same result than ours, however, we give here a different and slightly shorter proof.  These two results
allow to understand the evolution of Wigner measures of  families of solutions to
\aref{eq:schro} in Section~\ref{sec:measures}.

$ $

This paper is the first step in the understanding of  symmetric
conjoint crossings.  For applications in quantum chemistry, it would be interesting to be able to treat general potential $V(q)$ presenting a conjoint crossing of symmetric eigenvalues on a codimension~$2$ submanifold. Then,  it is likely that the potentials $V_{JT}$ and $V_{PJT}$ will play a crucial role as toys model. 
 This motivates the work performed in this paper. Besides, 
the application to quantum chemistry that we have in mind concerns the elaboration of algorithms 
modelizing the evolution of the Wigner transform of families of solutions of Schr\"odinger equation~\aref{eq:schro}. These algorithms have  been extensively developed since the 70's and the pioneer work of Tully and Preston (see \cite{TP}). The algorithms of  \cite{LT} and \cite{FL}  which are the first one presenting a mathematical proof  relies on normal forms results of~\cite{CdV1}, \cite{CdV2} and \cite{F07} and on the precise analysis of the propagation of Wigner measures through the crossing performed in \cite{FG03} and \cite{F03}. Therefore,   
the result of this article concerning
the branching of Wigner measures has important consequences for
numerics in quantum chemistry:  one can construct an
algorithm for the potential $V_{PJT}(q)$  exactly as in  \cite{LT} and its convergence is
guaranteed by the Theorem~\ref{theo:measures} in the same manner
than the algorithm of \cite{LT} relies on the branching 
theorem of \cite{FG02}.  The diagrams illustrating the statement of the main result in the next section are obtained
that way.

%%%%%%%%%%%%%%%%%%%%%%%%%%%%%%%%%%%%%%%%%%%%%%%%%%%%%%%%%%%%%%%%%%%%%%%%%%%%%%%%%%%%%%%%%%%%%%%%%%%%%

\section{Main result}

Let us  begin with some notations.  For
$\ell\in\{0,+1,-1\}$, we denote by $\Pi^\ell(q)$ the
eigenprojector associated with $\ell \mid q\mid$ and we set
$$\lambda^\ell(q,p,\tau)=\tau+{|p|^2\over 2}+\ell|q|.$$
The indices $\ell\in\{0,+1,-1\}$ will be  sometimes shortened into
$\pm$ whenever $\ell=\pm1$. We consider the Hamiltonian vector fields
$$H_{\lambda^\ell}=\partial_\tau+p\cdot\nabla _q -\ell{q\over|q|}\cdot\nabla_{p},\;\;\ell\in\{0,+1,-1\}
$$ and their integral curves $\rho^\ell_s$  which are called classical trajectories.
They are of the form  $$\rho^\ell_s=(q_s^\ell, s +t_0, p_s^\ell,\tau_0) \;\;{\rm  with}\;\; 
\dot q^{\ell}_s=p_s^{\ell},\;\;  \dot
p^{\ell}_s=\ell\;{q_s^\ell\over |q_s^{\ell}|}\;\;{\rm and}\;\;\rho^\ell_0=(q_0,t_0,p_0,\tau_0).$$
Because of the singularity in $q=0$, something has to be said for the $\pm$ mode: it is proved in
\cite{FG02} that under the assumption
\begin{equation}\label{generic}
p_0\not=0,
\end{equation}
 there exists a unique curve $\left(q^{\pm}_s,p^{\pm}_s\right)$
solution of the system of ordinary equations above and such that~$q_0^\pm=q_0$ and $p^\pm_0=p_0$. Moreover,  if $q_0=0$, the
$+$~trajectory  (resp. $-$~trajectory) issued from~$(q_0,p_0)$
smoothly continues the $-$~trajectory (resp. $+$~trajectory)
arriving at $(q_0,p_0)$. Finally, we emphasize a
 specific  feature of these trajectories:
there exists an hypersurface containing all the curves which pass
through the crossing set, namely the singular hypersurface
$$I=\{ q\wedge p=0\}$$
where for  two vectors $z=(z_1,z_2)$ and $z'=(z'_1,z'_2)$ of
$\R^2$, $z\wedge z'$ denotes their wedge product:
$$z\wedge z'=z_1z_2'-z_2z_1'.$$
Observe that $I$ is a  smooth hypersurface close to points $\rho_0=(q_0,t_0,p_0,\tau_0) $ such that $p_0\not=0$.
Note also that this assumption implies that the classical trajectories are transverse to the crossing set. It is under this assumption that most results on crossings are obtained (for example in  \cite{Ha94}, \cite{HJ}, \cite{FG02} or \cite{FG03}; for a result in case of tangency see \cite{DFJ}).  

\ni We work in space-time variables and call {\it crossing set} the points of the phase space where
there  is an eigenvalue crossing namely the set
$$S=\left\{\tau+{|p|^2\over
2}=0,\;\;q=0\right\}=\{
q=0\}\cap\Sigma,$$ where $\Sigma$ denotes the characteristic set:
$$\Sigma=\Sigma^+\cup\Sigma^-\cup \Sigma^0\;\;{\rm 
with}\;\;
\Sigma^\ell=\{\lambda^\ell(q,p,\tau)=0\}$$ being
the energy surface for the mode $\ell$. 
\ni The Wigner measures of families $(\psi^\eps)$ are
supported on the energy surfaces. More precisely, it is proved in
\cite{GMMP} that outside the crossing set, one has
$$\mu=\mu^+\,\Pi^+ + \mu^- \,\Pi^-+\mu^0\,\Pi^0$$
where for $\ell\in\{0,+1,-1\}$, the measures $\mu^\ell$ are scalar
positive Radon measure such that for all $a\in{\mathcal
C}_0^\infty(\R^6,\C)$,
$$\left(\op_\eps(a\Pi^\ell)\psi^\eps,\psi^\eps\right)_{L^2(\R^3,\C^3)}\td_\eps,0
\int_{\R^6}a(q,t,p,\tau)\, \d\mu^\ell(q,t,p,\tau).$$ Besides, the
measures $\mu^\ell$ are supported on $\Sigma^\ell$ and invariant
through the Hamiltonian flow associated with $H_{\lambda^\ell}$ for all $\ell\in\{0,-1,+1\}$.

\ni It is well known that the transitions between two modes cannot
be described only in terms of Wigner measures (see \cite{FG02}). The phenomenon of
the transitions is more intricate and requires a second level of
observation. Indeed,  the transitions are determined by the way
the families $(\psi^\eps(t))$ concentrate on the trajectories
arriving at the crossing with respect to the scale $\sqrt\eps$.
For this reason, one uses two-scale Wigner measures associated with the hypersurface $I$ because all the trajectories passing through $S$ are contained in $I$.

\ni Following \cite{Mi} and \cite{FG02}, we consider observables
$a=a(q,t,p,\tau,\eta)\in{\mathcal C}^\infty(\R^{7}, \C^{3,3})$ which are compactly
supported outside $\{p=0\}$ in the variables $(q,t,p,\tau)$, uniformly with respect
to the variable $\eta$ and coincide for $|\eta|>R_0$ for some
$R_0>0$ with a function homogeneous of degree $0$ in the variable
$\eta$. We will denote by ${\mathcal A}$ the vector space of such
functions. Then, when given the hypersurface $I$, one associates
with $a$ the pseudodifferential operator
$$\op_\eps^{I}(a)=\op_\eps\left(a\left(q,t,p,\tau,{q\wedge
p\over\sqrt\eps}\right)\right).$$ The change of variables
$$(q,p)\mapsto(\sqrt\eps q,\sqrt\eps p)$$
and Calder\'on-Vaillancourt Theorem show that this family of
operator is uniformly bounded with respect to $\eps$ on
$L^2_{loc}(\R^3)$. A two-scale Wigner measure associated
with the  concentration of the family $(\psi^\eps)$ on $I$ is a
positive Radon measure on the compactified normal bundle to $I$
such that
$$\displaylines{\qquad
\left(\op_\eps^I(a)\psi^\eps,\psi^\eps\right)_{L^2(\R^3,\C^3)}\td_\eps,0\hfill\cr\hfill
\int_{\overline N (I)}a(q,t,p,\tau,\eta)\d\nu(q,t,p,\tau,\eta)+ \int _{I^c} a\left(q,t,p,\tau,{ q\wedge p
\over |q\wedge p|}\,\infty\right) \d \mu(q,t,p,\tau).\qquad\cr}$$
We recall that the fibres of the normal bundle to $I$ is obtained above a point
$\rho$ of $I$ by quotienting the vector space
$T_\rho(T^*\R^3)$ by the tangent to $I$, $T_\rho I$. The fibre
above $I$ is a dimension $1$ vector space $N_\rho(I)$ and one
gets the fibre of the compactified bundle to $I$ by adding two
points at infinity; then, this fiber is isomorphic to
$\overline\R$. The function $a$ being homogeneous in the variable
$\eta$, it defines a function on $\overline N_\rho(I)$. We point
out that $\mu$ above $I$ is the projection of $\nu$ on $I$:
  $$\displaylines{  \mu(q,t,p,\tau){\bf
1}_I=\int_{\overline\R}\nu(q,t,p, \tau,\d\eta).\cr}$$
Since $\mu$ is determined outside $I$ by the Wigner measures of the data (because of the invariance of $\mu^\ell$ through $H_{\lambda^\ell}$ outside $S$), we focus on what happens above $I$ and thus on  the measure $\nu$.

\ni Similarly to Wigner measures, the two-scale Wigner measures
satisfy localization and propagation properties: it is proved in
\cite{FG02} that outside $S$,
$$\nu=\nu^0\,\Pi^0+\nu^+\,\Pi^+ +\nu^-\,\Pi^-$$
where for any $\ell\in\{0,+1,-1\}$, the scalar positive Radon
measures $\nu^\ell$ are supported on $\Sigma^\ell$ and propagate
along the linearized Hamiltonian flow induced on $\overline N (I)$ by the
Hamiltonian vector fields~$H_{\lambda^\ell}$. We are now able to describe the branching of two scale Wigner
measures close to a point $\rho_0=(q_0,t_0,p_0,\tau_0)$ such that
$p_0\not=0$. Indeed, since $p_0\not=0$, the
classical trajectories are transverse to the crossing set $S$ and
for any $\ell\in\{0,+1,-1\}$, the measures $\nu^\ell$ have traces
on $S$. Let us denote by $\nu^\ell_{in}$ (resp. $\nu^\ell_{out}$)
the traces of $\nu^\ell $ on the in-going (resp. out-going)  side
of $S$. For expressing the link between the traces, we need to
choose coordinates on $\overline N(I)$: for $\rho=(q,t,p,\tau)\in
I$, we characterize the class in $N_\rho(I)$
of the vector $\delta\rho=(\delta
q,\delta t,\delta p,\delta \tau)$ of $T_\rho I$ by the coordinate
$$\eta=\delta q\wedge  p .$$ These traces are  linked as
stated in the following theorem.

\begin{theo}\label{theo:measures}
Above $\rho_0=(q_0,t_0,p_0,\tau_0)\in S$ with $p_0\not=0$ and under the condition that $\nu^{0}_{in}$, $\nu^{+}_{in}$ and~$\nu^{-}_{in}$ are two by two mutually singular on $\{|\eta|<+\infty\}$, 
one has
\begin{equation}\label{link}
\left(\begin{array}{c}
\nu^+_{out}\\\nu^-_{out}\\\nu^0_{out}\end{array}\right)= 
\left(\begin{array}{ccc}
(1-T)^2 & T^2 & 2T(1-T) \\
T^2 & (1-T)^2 & 2T (1-T) \\
2T (1-T) & 2T(1-T) & (1-2T)^2\end{array}\right)
\left(\begin{array}{c}
\nu^+_{in}\\\nu^-_{in}\\\nu^0_{in}\end{array}\right)
\end{equation} where
\begin{equation}\label{trscoef}
T(p,\eta)={\rm exp}\left({-{\pi \eta^2\over 2|p|^3}}\right).
\end{equation}
\end{theo}

\ni The coefficient $T(p,\eta)$ is very close to the Landau-Zener coefficient ${\rm exp}\left({-{\pi \eta^2\over |p|^3}}\right)$ which appears for crossings of two modes (see \cite{La}, \cite{Ze}, \cite{HJ} and \cite{FG02}). 

\ni We could also have expressed the transitions through the crossing set by considering the two-scale Wigner measures associated with the concentration of $(\psi^\eps)$ on the sets $J^{\ell}_{in}$ (resp. $J^\ell_{out}$) consisting of all the trajectories for the mode $\ell$ arriving to (resp. going out of) the crossing set for $\ell\in\{0,+1,-1\}$. These measures have traces on $S$ that can be identified to the measures $\nu^\ell_{in/out}$. The link is then given by \aref{link}. This point of view  (which is the one of \cite{FG02})  is more intrinsic and applies in situations where there is no trivial sets containing all the trajectories passing through the crossing. This is a real issue for general symmetric conjoint crossing (see \cite{FR}).

\ni We have reduced the dimension space to $2$ but we could have assumed the dimension to be greater than $2$ with a $d$-dimensionnal  Laplacian and a  potential $V_{PJT}(q_1,q_2)  $ depending only on $q_1$ and $q_2$. Our result extends to this situation replacing $I$ by $\{(p_1,p_2)\wedge(q_1,q_2)=0\}$ and $p$ by $(p_1,p_2) $ in the non-degeneracy condition~\aref{generic} and in the transition coefficient~\aref{trscoef}. It also extends to the situation where the Laplacian is replaced by a Fourier multiplier $A(D)$. One then need to turn $p$ into $\nabla A(p)$ in   \aref{generic} and~\aref{trscoef}. In that situation, the existence of a submanifold satisfying similar properties than $I$ is non trivial.

\ni  The proof of this theorem relies on a normal form result which is
the subject of Section~\ref{sec:nf}. One reduces the analysis of
the crossing to a scattering problem which is studied in
Section~\ref{sec:scat}. Provided these two results, the proof of
Theorem~\ref{theo:measures} is performed in the last section following the method introduced in \cite{FG02}
(slightly modified to take into account the difficulties induced by the presence of  the third mode).

$ $

\ni Before passing to the proof of the Theorem, we illustrate our result by some numerics.  We choose a Gaussian  initial data  supported on the $+$ level of the form 
$$\psi^\eps_0(q)=  (\eps\pi)^{-1/2}
\exp\!\left(- {\textstyle\frac{1}{2\eps}}|q-q_0^\eps|^2 +
{\textstyle\frac{i}{\eps}}\, p_0\cdot(q-q_0^\eps) \right) e_+(q)$$
where $$e_+(q)\in{\rm Ran}\Pi^+(q),\;\;\eps=10^{-2},\;\;p_0=(-1,0),\;\;q_0^\eps=\sqrt\eps\left(5, {1\over 2}\right).$$
The time evolves between $t=0$ and $t= {\pi\over 4}$ and we  calculate  an approximate value  of the  population on each mode: 
$$n^\ell(t)=\int _{\R^2} |\Pi^\ell(q)\psi^\eps(q,t)|^2 dq,\;\;\ell\in\{-1,0,+1\}.$$
Following \cite{LT} and \cite{FL}, we perform  a surface hopping algorithm:
\begin{itemize}
\item we sample the Wigner function, 
\item we propagate  weighted points $(q^+_s,p^+_s, w)$ along the classical trajectory until the distance to the gap $\{q=0\}$ is minimal, 
\item at that point $(q^+_{s^*},p^+_{s^*},w)$  one generates new trajectories for each mode transporting part of the  weight~$w$:
$$(q^+_s,p^+_s,(1-T^*)^2w),\;\; (q^-_s,p^-_s,{T^*}^2w),\;\;(q^0_s,p^0_s,2T^*(1-T^*)w),\;\;s>s^*$$
where the result of Theorem~\ref{theo:measures} motivates the choice of the  transition coefficient 
$$T^*=\exp\left(-{\pi|q^0_{s^*}\wedge p^0_{s^*}|^2\over2\,\eps|p^0_{s^*}|^3 }\right).$$
\end{itemize}
Arguing as in  \cite{LT}, one can prove the convergence of the algorithm. No estimation on the rate is available for the moment since it requires an  improvement of the normal form result which is out of reach for the moment.   In order to illustrate numerically the validity of the algorithm, we calculate a reference solution based on a grid discretization and  a Strang splitting scheme.  This scheme is known to have a very fast rate of convergence (see \cite{LST}, \cite{FL} and \cite{FL2}).
Then, we compare the two outputs in the figure below: the initial population  which was on the $+$ level splits  after passing once through the crossing.

\begin{figure}[h]
\begin{center}
\scalebox{0.4}{\includegraphics{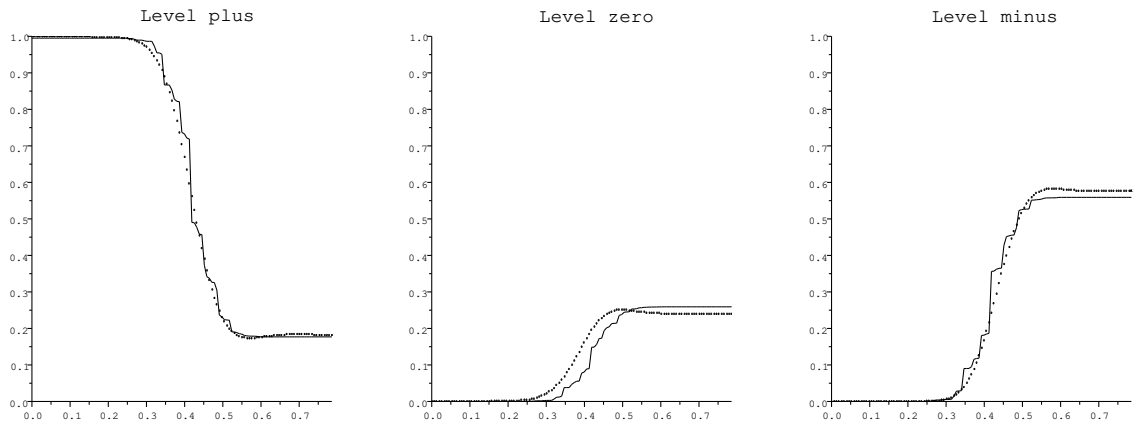}}
\caption{Evolution of the population on each level for surface hopping algorithm (full line) and for Strang splitting scheme (dotted line).}
\end{center}
\end{figure}

%%%%%%%%%%%%%%%%%%%%%%%%%%%%%%%%%%%%%%%%%%%%%%%%%%%%%%%%%%%%%%%%%%%%%%%%%%%%%%%%%%%%%%%%%%%%%%%%%%%%%%%%%%%%

\section{Normal form}\label{sec:nf}

In this section, we first prove a microlocal normal form which holds in the phase space close to a point of the crossing set. Then, in the second subsection, we write a two-microlocal normal form which is only valid in neighborhoods of size $\sqrt\eps $ of the hypersurface $I$.

\subsection{The microlocal normal form}

The phase  space $\R^6$ has the structure of the cotangent
space $T^*\R^3$. It is endowed with a symplectic structure
given by the $2$-form
$$\omega=\d\gamma=\d p\wedge\d q+\d\tau\wedge\d t$$
 where $\gamma$ is the Liouville $1$-form $\gamma=p\,\d q+\tau\d t$.
 A local change of sympletic coordinates $\kappa$ is called {\it  canonical transform}
 and
 one can  associate with it a unitary bounded operator of
 $L^2(\R^3)$ compatible with the change of coordinates. Such
 an operator $K$ is called {\it Fourier integral operator} associated
 with $\kappa$. It satisfies Egorov's Theorem
\begin{equation}\label{eq:Egorov}\forall
a\in{\cal C}^\infty_0(\R^6),\;\;
K^*\op_\eps\left(a\right)K=\op_\eps(a\circ\kappa)+O(\eps^2)\;\;{\rm
in}\;\;{\cal L}(L^2(\R^3)) . \end{equation} \ni   We point out that the remainder is better than expected because we have chosen to use Weyl quantization.  The reader
will find in~\cite{Rb} a complete analysis of Fourier integral
operators.

$ $
 
\ni  We consider $\Omega$ a neighborhood of $\rho_0=(q_0,t_0,p_0,\tau_0)\in S$ with $p_0\not=0$ and we suppose $\Omega$ small enough so that $p\not=0$ in $\Omega$.  Then, for $\ell\in \{0,+1,-1\}$,  we denote by $J^\ell_{in}$ (resp. $J^\ell_{out}$) the sets consisting on all the classical trajectories for the mode $\ell$ entering in (resp. ougoing of) $S\cap \Omega$. We point out that by the property of continuation of the trajectories, the sets $\left(J^\pm_{in}\cup J^\mp_{out}\right)\cap \Omega$ are submanifolds of $T^*\R^3$ just as $\left(J^0_{in}\cup J^0_{out}\right)\cap \Omega$ does.
 Our result is the
following

\begin{theo}\label{theo:nf}
Let $\rho_0=(q_0,t_0,p_0,\tau_0)\in S$ be such that $
p_0\not=0$.  Then, there exists a local canonical transform
$\kappa$ from a neighborhood of $\rho_{0}$ into some neighborhood
$\tilde\Omega$ of $0$,
$$\kappa: (q,t,p,\tau)\mapsto(s,z,\sigma,\zeta),\;\;\kappa(\rho_0)=0,$$
a Fourier integral operator $K$
associated with $\kappa$ and an invertible  matrix-valued function
$A$ such that if $\psi^\eps$ is a family of solutions to
\emph{\aref{eq:schro}} for some initial data $\psi^\eps_0$ uniformly
bounded in~$L^2(\R^2)$, then
$$v^\eps=K^*\op_\eps\left(A\right)\psi^\eps$$
satisfies for all $\phi\in{\cal C}_0^\infty(\tilde\Omega)$,
\begin{equation}\label{systreduit}
\op_\eps(\phi)\,\op_\eps\Bigl(-\sigma+V_{PJT}\bigl(s,\alpha(\zeta,\sigma)\,
z_1\bigr)\Bigr)\,v^\eps= O(\eps)
\end{equation}
\noindent  in $L^2(\R ^{3})$ where  $\alpha$ is smooth and does not vanish in $\tilde\Omega$.

\noindent Moreover, 
\begin{eqnarray}
\label{eq:newI} & \kappa(I)=\{z_1=0\},&\\
\label{eq:newJ+p}
& \kappa(J^\pm_{in})=\{\sigma\pm s=0,\;\; z_1=0,\;\;s\leq 0\},&\\
\label{eq:newJ+f} &\kappa(J^\pm_{out})=\{\sigma\mp s=0,\;\;
z_1=0,\;\;s\geq 0\}.&\end{eqnarray}
\end{theo}

\begin{dem}
The proof of this Theorem proceeds in two steps. First, we find the canonical transform~$\kappa$. Then, we use the properties of the matrix $V_{PJT}$ to construct $A$ such that 
$$\left[A\left(\left(\tau +{|p|^2\over 2}\right){\rm Id} +V_{PJT}(q_1,q_2)\right)A^*\right]\circ \kappa^{-1}=-\sigma\, {\rm Id} +V_{PJT}(s,z).$$
Then, we obtain \aref{systreduit} by the Egorov theorem \aref{eq:Egorov}. 
 
$ $

\ni{\bf First step: the canonical transform}.  This first step relies on the analysis of the geometry of the crossing. 
We crucially use Proposition 6, p.148 in \cite{FG02}.

\begin{prop}\label{prop:cantrsf}
There exist a local canonical transform
$$\kappa: (q,t,p,\tau)\mapsto (s,z,\sigma,\zeta)$$
and non-vanishing smooth functions $\lambda$ and $\mu$ such that 
\begin{equation}\label{eq:coordinates}
\left\{ \begin{array}{l}
\sigma=\lambda(p,\tau)\left(\tau +{|p|^2\over 2}\right),\\
s= \lambda (p,\tau) {p\over |p |}\cdot q,\\
z_1=\lambda(p,\tau)\mu(p,\tau){p\over|p|}\wedge q.
\end{array}\right.\end{equation}
Besides, one can choose
 $\zeta=\zeta(p,\tau)$ such that $(p,\tau)\mapsto(\zeta,\sigma)$ is a diffeomorphism and 
 \begin{equation}\label{lambda}
 \lambda(p,\tau)_{|S}=|p|^{-1/2}.
 \end{equation}
\end{prop}

\noindent Let us shortly describe the proof of Proposition 6 in \cite{FG02}. Observing that $$\left\{\tau +{|p|^2\over 2} ,{p\over|p|}\cdot q\right\}=|p|,$$
we obtain the existence of $\lambda$ such that  $s$ and $\sigma$ satisfy the bracket condition $\{\sigma,s\}=1$. Therefore we necessarily have \aref{lambda}. Besides, for any function $\mu=\mu(p,\tau)$, one has 
$$\left\{ \lambda(p,\tau)\left(\tau +{|p|^2\over 2}\right), \lambda(p,\tau)\mu(p,\tau) {p\over|p|}\wedge q\right\}=0.$$
 Therefore, it is enough to find  $\mu$ such that 
 $$\left\{\lambda(p,\tau) {p\over|p|}\cdot q\;,\;\lambda(p,\tau)\mu(p,\tau) {p\over|p|}\wedge q\right\}=0.$$
  Once $\mu$ is built,  one completes the symplectic system $(\sigma,s,z_1)$ by Darboux Theorem.

  $ $

\noindent By the definition of $S$ and $I$, one has $\kappa(S)=\{\sigma=s=z_1=0\}$ and $\kappa(I)=\{z_1=0\}$. Moreover, 
since $\lambda\not=0$ (since $p\not=0$ in $\Omega$), we have  $\kappa(\Sigma)=\{\sigma(\sigma^2-s^2-z_1^2)=0\}$.  In order to find the precise equations of $\kappa(J^\pm_{in})$ and $\kappa(J^{\pm}_{out})$, we observe that 
\begin{equation}\label{eq:Ietcie}
\kappa(J^\pm_{in/out})\subset
\{\sigma^2=s^2,\;\;z_1=0\}
\end{equation}
and  we consider the Hamiltonian vector fields  
$H_{\lambda^\pm}$. We have outside $S$
$$H_{\lambda^\pm}=\partial_t+p\cdot \nabla_q \mp\, {q\over |q|}\cdot\nabla_{p}.$$
Observing that if the trajectory $(q^\pm_s,p^\pm_s)$ passes through $S$ at $s=0$ we have 
$$q^\pm_s =s\,p_0 +o(s)$$
we deduce
$$H_{\lambda^\pm}(\rho^\pm_s)=\partial_t +p\cdot \nabla _q \mp\,{\rm sgn}(s) {p_0\over|p_0|}\cdot\nabla_{}+o(1).$$
By the definition of $\sigma$ and $s$, we have 
$$H_{\sigma\pm s}=\lambda\left[\partial_t + p\cdot \nabla_q  \pm \,{p\over|p|}\cdot\nabla _{}\right]\;\;{\rm on}\;\; S.$$
Therefore, we obtain 
$$H_{\lambda^+}(\rho_s^+)\td_{s},{0^-} \partial_t +p_0\cdot \nabla_q + {p_0\over|p_0|}\cdot \nabla_{±q}=\lambda^{-1} H_{\sigma+s} .$$
This equation implies that $\kappa(J^{+}_{in})\subset\{ \sigma+s=0\}$, whence the equation of $\kappa(J^+_{in})$ by \aref{eq:Ietcie} and dimension considerations. One argues similarly for $J^+_{out}$, $J^{-}_{in}$ and $J^-_{out}$.

$ $

\ni{\bf Second step: the gauge transform}. We now transform the symbol  $\left(\tau +{|p|^2\over 2}\right){\rm Id} +V(q)$ in order to see the coordinates $(s,z,\sigma,\zeta)$ in the symbol.
We first use the matrix 
$$M={1\over \sqrt 2}\left(\begin{array}{ccc}
1 & 1 & 0 \\ 1 & -1 & 0 \\ 0 & 0 & \sqrt 2 \end{array}\right)$$
and observe that $M=M^*=M^{-1}$ and $MV(q)M=W(q)$ with
$$W(q)=\left(\begin{array}{ccc} 0 & q_1 & q_2 \\ q_1 & 0 & 0 \\ q_2 & 0 & 0 \end{array}\right).$$
We then use the following lemma which comes from a straightforward computation.

\begin{lem} 
 Consider the matrix $\displaystyle{ R(p)={1\over|p|}\left(\begin{array}{ccc} |p| & 0 & 0 \\ 0 & -p_1 & -p_2 \\ 0 & p_2 & -p_1\end{array}\right)}$ then  $R(p)$ is invertible for $p\not=0$ and we have
 $$ R(p) W(q)R(p)^* =- W\left( {p\over|p|}\cdot q,  {p\over|p|}\wedge q\right).$$
\end{lem}

$ $

We can now conclude the proof of Theorem~\ref{theo:nf}.
 Set 
$$P(q,p,\tau)=\left(\tau+{|p|^2\over 2}\right){\rm Id} +V(q)\;\;{\rm and}\;\;\tilde R(p)=M  R(p) M.$$
We have
$$\tilde R(p)P(q,p,\tau) \tilde R(p)^*=\left(\tau +{|p|^2\over 2} \right){\rm Id}+  V\left( {p\over|p|}\cdot  q ,  {p\over|p|}\wedge q\right).$$
Therefore, using the canonical transform of Proposition~\ref{prop:cantrsf}, we have
$$\left[\lambda(p,\tau) \tilde R(p)P(q,p,\tau)\tilde R(p)^*\right] \circ \kappa^{-1}(s,z,\sigma,\zeta)= \sigma\,{\rm Id}-V(s,\alpha(p,\tau) z_1).$$
Let us consider now a Fourier integral operator $U$  associated  with $\kappa$ and set
$$A(p,\tau)=\lambda(p,\tau) ^{1/2}\tilde R(p).$$
We have
$$ U \op_\eps(A)\op_\eps(P)\op_\eps(A^*)U^*= \op_\eps\left(\sigma\,{\rm Id} -V(s,\alpha(\sigma,\zeta) z_1)\right) +O(\eps)$$
in ${\mathcal L}\left(L^2(\R^3)\right)$.
\end{dem}

\subsection{The 2-microlocal normal form}\label{sec:2micro}

In this section, we want to ameliorate the normal form in open sets localized at a distance of order $\sqrt \eps $ of the hypersurface $I$.  For simplicity, we set 
\begin{equation}\label{eq:L}
J=\left(\begin{array}{ccc}
1 & 0 & 0 \\ 0 & -1 & 0 \\ 0 & 0 & 0 \end{array}\right),\;\; L= {1\over \sqrt 2}\left(\begin{array}{ccc}
0 & 0 & 1 \\ 0 & 0 & 1 \\ 1 & 1 & 0 \end{array}\right).
\end{equation}
We shall take advantage of the fact that  if $a\in{\mathcal C}_0^\infty(\R^{7})$, 
$$\op_\eps^I(a) \op_\eps \left( -\sigma {\rm Id} +V_{PJT}(s,\alpha(\sigma,\zeta) z_1)\right)=\op_\eps^I (a) \op_\eps^I \left( -\sigma {\rm Id} +V_{PJT}(s,\alpha(\sigma,\zeta)\sqrt\eps \eta)\right).$$
We set
\begin{eqnarray}\label{def:Q}
Q & = & -\sigma {\rm Id} +s J +\sqrt \eps\,\alpha(\sigma, \zeta) \eta L,\\
\label{def:Q0}
Q_0 &  =  &  -\sigma {\rm Id} +s J +\sqrt \eps \,\alpha(0, \zeta) \eta L.
\end{eqnarray}
We prove the following result.

\begin{prop}\label{prop:nf2} For every ball $B$ of $\R_\eta$,
there exist smooth matrix-valued functions $C$  and $\tilde C$  such that   for all $a\in{\mathcal C}_0^\infty (\R^6\times B)$,
$$\left\|\op_\eps^I(a)\left[ \op_\eps^I({\rm Id}+\sqrt\eps\, C) \op_\eps(Q) -\op_\eps(Q_0) \op_\eps^I({\rm Id}+\sqrt\eps\,\tilde C) \right] \right\|_{{\mathcal L}(L^2)}=O(\eps).$$
\end{prop}

\begin{dem} We want to realize 
$$\op_\eps^I({\rm Id}+\sqrt\eps C) \op_\eps(Q) =\op_\eps(Q_0) \op_\eps^I({\rm Id}+\sqrt\eps\,\tilde C) +O(\eps).$$
For this, we successively consider the terms of order $1$ and $\sqrt \eps$ in the development of the preceeding equality by symbolic calculus. The terms of order $0$ are equal and the equality for the terms of 
 order $\sqrt \eps$ reads
$$\alpha(\sigma,\zeta) \eta \, L-\alpha(0,\zeta) \eta \, L +\sigma (\tilde C -C) +s(C J-J\tilde C)=0.$$
We choose $\tilde C=-\,^tC$ with 
$$C=f \left(\begin{array} {ccc} 0 & 0 & 1 \\ 0 & 0 & 1 \\ 0 & 0 & 0 \end{array}\right)$$
so that the second equation writes in view of $CJ=0$  and $J\tilde C=- \,^t(CJ)=0$
$$\left[\eta (\alpha(\sigma,\zeta)-\alpha(0,\zeta))-\sqrt 2\, \sigma f\right]L=0.$$
Hence
$$
f(\sigma,z,\eta)  = { \eta\over\sqrt 2}{\alpha(\sigma,\zeta) -\alpha(0,\zeta)\over\sigma}.$$
This closes the proof of the Proposition.

\end{dem}
%%%%%%%%%%%%%%%%%%%%%%%%%%%%%%%%%%%%%%%%%%%%%%%%%%%%%%%%%%%%%%%%%%%%%%%%%%%%

\section{Scattering matrix for the pseudo Jahn-Teller model}\label{sec:scat}

In this section, we analyze the  differential system
\begin{equation} \label{PJTsystem}
-i\partial_s\left(\begin{array}{c} v_+ \\ v_- \\ v_0 \end{array}\right)=\left(\begin{array}{ccc} s & 0 & z/\sqrt{2} \\ 0 & -s & \bar{z}/\sqrt{2} \\ \bar{z}/\sqrt{2} & z/\sqrt{2} & 0 \end{array}\right)\left(\begin{array}{c} v_+ \\ v_- \\ v_0 \end{array}\right)=V_{PJT}(s,z)\left(\begin{array}{c} v_+ \\ v_- \\ v_0 \end{array}\right)
\end{equation}
where $z$ is a non-zero complex parameter. We prove the following scattering result.

\begin{theo}\label{th:scat}
There exists $\left(\alpha^-,\alpha^+,\alpha^0\right)$ and
$\left(\omega^-,\omega^+,\omega^0\right)$  such that
$$\displaylines{
v_\ell(s)={\rm exp}\left[\ell\,i \left({s^2\over 2}+{|z|^2\over
2}{\rm ln}|s|\right)\right] \alpha^\ell +O\left(\frac{|z|^2}{s}\right)\;\;{\rm as}\;\;
s\rightarrow -\infty,\cr v_\ell(s)={\rm exp}\left[\ell\,i
\left({s^2\over 2}+{|z|^2\over 2}{\rm ln}|s|\right)\right]
\omega^\ell +O\left(\frac{|z|^2}{s}\right)\;\;{\rm as}\;\; s\rightarrow +\infty.\cr}
$$
Besides
$$\left(\begin{array}{c} \omega^+\\ \omega^-\\
\omega^0\end{array}\right)=S(z)
\left(\begin{array}{c} \alpha^+\\ \alpha^-\\
\alpha^0\end{array}\right)$$ with
$$S(z)=\left(\begin{array}{ccc}
e^{-\pi|z|^2/2} & i\Omega(z)^2[1-e^{-\pi|z|^2/2}] & \sqrt{2}e^{i\pi/4}\Omega(z)\theta(z) \\
-i\overline{\Omega(z)}^2[1-e^{-\pi|z|^2/2}] & e^{-\pi|z|^2/2} & -\sqrt{2}e^{-i\pi/4}\overline{\Omega(z)}\theta(z) \\
-\sqrt{2}e^{-i\pi/4}\overline{\Omega(z)}\theta(z) & \sqrt{2}e^{i\pi/4}\Omega(z)\theta(z) & 2e^{-\pi|z|^2/2}-1
\end{array}\right)$$
and
$$
\Omega(z)=\frac{z}{|z|}
\frac{\Gamma\left(1-i\frac{z^2}{4}\right)}{\left|\Gamma\left(1-i\frac{z^2}{4}\right)\right|},
\quad \theta(z)=\sqrt{1-e^{-\pi z^2/2}}e^{-\pi z^2/4}.
$$
\end{theo}

Note that the indices $+$, $-$ and $0$ are not directly connected with the modes $\lambda^\ell$. However their connection will be clarified in Section \ref{subsec:finitedist}.

\subsection{Providing two independent solutions with contour integrals}

Classical theorems give that the solutions form a tridimensional linear space of analytic functions. The parity block decomposition of $V_{PJT}(s,z)$
induces that the space of solutions split into two subspaces: a bidimensional one with $v_\pm$ even and $v_0$ odd, and a unidimensional one with $v_\pm$ odd and $v_0$ even. For an even (respectively odd) analytic function $v$, there exists an analytic function $f$ such that $v(s)=f(s^2/2)$ (respectively $v(s)=sf(s^2/2)$). Hence the system, in the case of the bidimensional subspace,
\begin{equation} \label{system1}
-i\left(\begin{array}{ccc} \partial_t & 0 & 0 \\ 0 & \partial_t & 0 \\ 0 & 0 & 1+2t\partial_t \end{array}\right)\left(\begin{array}{c} f_+ \\ f_- \\ f_0 \end{array}\right)=V_{PJT}(1,z)\left(\begin{array}{c} f_+ \\ f_- \\ f_0 \end{array}\right).
\end{equation}

\ni To solve \aref{system1}, we will use the method of contour integral (see for instance \cite{[Dieu]} Chap. XV): we look for solution of the form
\begin{equation} \label{contourintegral}
\left(\begin{array}{c} f_+(t) \\ f_-(t) \\ f_0(t) \end{array}\right)=\int_\Lambda e^{it\tau}\left(\begin{array}{c} F_+(\tau) \\ F_-(\tau) \\ F_0(\tau) \end{array}\right)d\tau
\end{equation}
where $\Lambda$ is an oriented path of $\C$ and the $F_*$'s are functions holomorphic in a neighbourhood of $\Lambda$. Injecting \aref{contourintegral} into \aref{system1} we get
$$
-\int_\Lambda e^{it\tau}i\left(\begin{array}{ccc} i\tau & 0 & 0 \\ 0 & i\tau & 0 \\ 0 & 0 & 1+2it\tau \end{array}\right)\left(\begin{array}{c} F_+(\tau) \\ F_-(\tau) \\ F_0(\tau) \end{array}\right)d\tau=\int_\Lambda e^{it\tau}V_{PJT}(1,z)\left(\begin{array}{c} F_+(\tau) \\ F_-(\tau) \\ F_0(\tau) \end{array}\right)d\tau
$$
and after integrating by parts
$$
\displaylines{ \left[e^{it\tau}\left(\begin{array}{c} 0 \\ 0 \\
-2i\tau F_0(\tau) \end{array}\right)\right]_\Lambda+\int_\Lambda
e^{it\tau}\left(\begin{array}{c} 0 \\ 0 \\ 2i[\tau
F'_0(\tau)+F_0(\tau)]
\end{array}\right)d\tau\hfill\cr\hfill=\int_\Lambda
e^{it\tau}\left(\begin{array}{ccc} -\tau+1 & 0 & \frac{z}{\sqrt{2}}
\\ 0 & -\tau-1 & \frac{\bar{z}}{\sqrt{2}} \\
\frac{\bar{z}}{\sqrt{2}} & \frac{z}{\sqrt{2}} & i
\end{array}\right)\left(\begin{array}{c} F_+(\tau) \\ F_-(\tau) \\
F_0(\tau) \end{array}\right)d\tau.\cr}
$$
Hence, it is enough to choose the $F_*$'s such that
$$
\left(\begin{array}{c} 0 \\ 0 \\ 2i\tau F'_0(\tau) \end{array}\right)=\left(\begin{array}{ccc} -\tau+1 & 0 & \frac{z}{\sqrt{2}} \\ 0 & -\tau-1 & \frac{\bar{z}}{\sqrt{2}} \\ \frac{\bar{z}}{\sqrt{2}} & \frac{z}{\sqrt{2}} & -i \end{array}\right)\left(\begin{array}{c} F_+(\tau) \\ F_-(\tau) \\ F_0(\tau) \end{array}\right)
$$
and $\Lambda$ such that
\begin{equation} \label{contour}
\left[e^{it\tau}\tau F_0(\tau)\right]_\Lambda=0.
\end{equation}

\begin{figure}[h]
\begin{center}
\scalebox{0.35}{\includegraphics{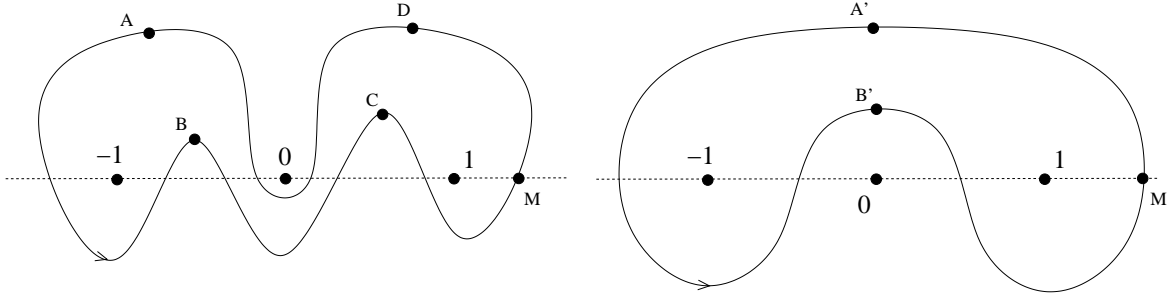}}
\caption{decomposition of the $\Lambda^\cup$ and $\Lambda^\cap$ contours \label{contour1}}
\end{center}
\end{figure}

\ni The system is exactly solvable and gives, up to a global constant,
$$
F_0(\tau)=\frac{e^{i\frac{|z|^2}{4}\left[\ln(\tau+1)-\ln(\tau-1)\right]}}{\sqrt{\tau}}, \quad
F_+(\tau)=\frac{z}{\sqrt{2}}\frac{F_0(\tau)}{\tau-1}, \quad
F_-(\tau)=\frac{\bar{z}}{\sqrt{2}}\frac{F_0(\tau)}{\tau+1}.
$$
Now $F_*(\tau)$ has three singularities, namely $\pm1$, $0$ and $F_*(\tau)$ is multivalued. More explicitly, considering a closed contour, the value of $F_*(\tau)$ is changed according to the following rule: if the contour goes $n_\pm$ times anticlockwise around $\pm1$ then $\ln(\tau\mp1)$ is changed into $\ln(\tau\mp1)+i2\pi n_\pm$ and if it goes $n_0$ times anticlockwise around $0$ then $\sqrt{\tau}$ is changed into $(-1)^{n_0}\sqrt{\tau}$. If we choose a closed contour $\Lambda$ that avoids the singularities of $F_*$ and solves \aref{contour}, the corresponding function \aref{contourintegral} is then well defined and is a solution of \aref{system1}, the only thing that one has to check is that it does provide a non trivial solution. We will focus on the following contours: they enclose only the singularities $\pm1$ and go once around each of them. There are exacly two of them depending on whether we stay below or above the singularity $0$ (in the first case the contour has the shape of a $\cup$ and the second of a $\cap$). The idea is now to use residue arguments to get explicit formulas for the contribution of the different singularities. Multiplying by $-\sqrt{2}/|z|$, we define our two solutions to be
\begin{equation} \label{solcontint}
\left(\begin{array}{c} v_+(s) \\ v_-(s) \\ v_0(s) \end{array}\right)=\int_\Lambda e^{i\frac{s^2}{2}\tau}e^{i\frac{|z|^2}{4}\left[\ln(\tau+1)-\ln(\tau-1)\right]}\left(\begin{array}{c} \frac{z}{|z|}\frac{1}{\tau-1} \\ \frac{\bar{z}}{|z|}\frac{1}{\tau+1} \\ \frac{\sqrt{2}}{|z|}s \end{array}\right)\frac{d\tau}{\sqrt{\tau}}
\end{equation}
where the $\Lambda^\cup$ and $\Lambda^\cap$ contours are decomposed according to Figure \ref{contour1} and the determination of the two logarithms and the square root is such that it coincides with the usual definition (for positive real numbers) at the point $M$.

\ni For the sake of clarity, the points $A$, $B$, $C$ and $D$ (respectively $A'$ and $B'$) have been split but they have to be thought of as the same geometric point. To study separately the different pieces, we will first show the following residue results.

\begin{lem}
Let $\beta$ be a non zero real number, $w(\tau)$ an holomorphic function of $\tau$ in a convex open neighborhood $\Omega$ of a point $\alpha$ of the real axis and $\Lambda$ an oriented path of $\Omega$ with both ends in the upper half plane. Then, in the limit $|s|\to+\infty$
\begin{enumerate}
\item
if $\Lambda$ does not surround $\alpha$,
$$
\int_\Lambda e^{i\frac{s^2}{2}\tau}e^{i\beta\ln(\tau-\alpha)}w(\tau)d\tau=O\left(e^{-\frac{s^2}{2}\delta}\right), \qquad \int_\Lambda e^{i\frac{s^2}{2}\tau}w(\tau)\frac{d\tau}{\sqrt{\tau-\alpha}}=O\left(e^{-\frac{s^2}{2}\delta}\right)
$$
where $\delta$ is the minimum of the imaginary part of the two endpoints ;
\begin{figure}[h]
\begin{center}
\scalebox{0.3}{\includegraphics{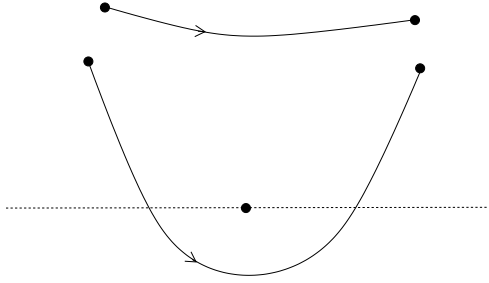}}
\caption{non-surrounding and surrounding contours\label{contourlemma}}
\end{center}
\end{figure}
\item
if $\Lambda$ does surround $\alpha$ once and anticlockwise,
\begin{equation}\label{logarithm}
\int_\Lambda
e^{i\frac{s^2}{2}\tau}e^{i\beta\ln(\tau-\alpha)}w(\tau)d\tau=-4ie^{i\beta\ln
2}e^{\pi\beta/2}\sinh(\pi\beta)\Gamma(1+i\beta)w(\alpha)\frac{e^{i\left[\frac{s^2}{2}\alpha-2\beta\ln|s|\right]}}{s^2}
\left[1+O\left(\frac{1}{s^2}\right)\right]
\end{equation}
\begin{equation} \label{squareroot}
\int_\Lambda e^{i\frac{s^2}{2}\tau}w(\tau)\frac{d\tau}{\sqrt{\tau-\alpha}}=2\sqrt{2\pi}e^{i\pi/4}w(\alpha)\frac{e^{i\frac{s^2}{2}\alpha}}{|s|}\left[1+O\left(\frac{1}{s^2}\right)\right]
\end{equation}
where the determination of the logarithm and of the square root is such that it coincides with the usual definition for positive real numbers. Furthermore \emph{\aref{logarithm}} is locally uniform in $\beta$.
\end{enumerate}
\end{lem}

\begin{dem}
\begin{enumerate}
\item
It follows from a crude estimate of the term $e^{i\frac{s^2}{2}\tau}$ for the straight line path.
\item
We deform the contour into the one presented in Figure \ref{contourlemma2} where the two vertical lines have been split for a better understanding {\it i.e.} $B$ and $E$ (respectively $C$ and $D$) are in fact the same geometric point. The radius of the small circle is chosen to be $\lambda(s)/s^2$ where $\lambda(s)\to0$ as $|s|\to+\infty$.
\begin{figure}[h]
\begin{center}
\scalebox{0.5}{\includegraphics{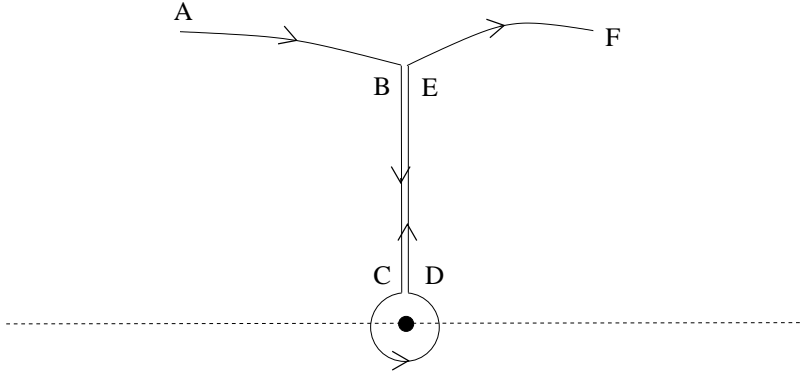}}
\caption{deformation of the surrounding contour\label{contourlemma2}}
\end{center}
\end{figure}

The pieces $AB$ and $EF$ are exponentially small by the first point. The piece $BC$ is $-e^{2\pi\beta}$ times (respectively equal to) the piece $DE$ for the logarithm (respectively the square root). Let us compute the piece $CD$ parameterizing the circle by $\alpha+\frac{\lambda(s)}{s^2}e^{i\theta}$ for $-3\pi/2\leq\theta\leq\pi/2$:
$$\displaylines{\qquad 
\int_C^D e^{i\frac{s^2}{2}\tau}e^{i\beta\ln(\tau-\alpha)}w(\tau)d\tau \hfill\cr\hfill=  \frac{\lambda(s)}{s^2}e^{i\frac{s^2}{2}\alpha}
e^{-2i\beta\ln|s|}e^{i\beta\ln\lambda(s)}
\int_{-3\pi/2}^{\pi/2}e^{i\frac{\lambda(s)}{2}e^{i\theta}}e^{-\beta\theta}w\left(\alpha+\frac{\lambda(s)}{s^2}e^{i\theta}\right)e^{i\theta}id\theta\cr\hfill
 =  O\left(\frac{\lambda(s)}{s^2}\right),\qquad\cr}$$
\begin{eqnarray*}
\int_C^D e^{i\frac{s^2}{2}\tau}w(\tau)\frac{d\tau}{\sqrt{\tau-\alpha}} & = & \frac{\sqrt{\lambda(s)}}{|s|}e^{i\frac{s^2}{2}\alpha}\int_{-3\pi/2}^{\pi/2}e^{i\frac{\lambda(s)}{2}e^{i\theta}}w\left(\alpha+\frac{\lambda(s)}{s^2}e^{i\theta}\right)e^{i\theta/2}id\theta \\
 & = & O\left(\frac{\sqrt{\lambda(s)}}{|s|}\right).
\end{eqnarray*}
Finally, we compute the contribution of the piece $DE$ parameterizing the vertical segment by $\alpha+2\frac{i}{s^2}u$ for $\lambda(s)/2\leq u\leq s^2\delta/2$:
$$\displaylines{
\int_D^E e^{i\frac{s^2}{2}\tau}e^{i\beta\ln(\tau-\alpha)}w(\tau)d\tau =
\frac{e^{i\frac{s^2}{2}\alpha}e^{-2i\beta\ln|s|}}{s^2}e^{i\beta\ln2}e^{-\pi\beta/2}2i
\int_{\lambda(s)/2}^{s^2\delta/2}e^{-u}u^{i\beta}w\left(\alpha+2\frac{i}{s^2}u\right)du \cr\hfill
  =  \frac{e^{i\frac{s^2}{2}\alpha}e^{-2i\beta\ln|s|}}{s^2}e^{i\beta\ln2}e^{-\pi\beta/2}2i\\
 \left[w(\alpha)\int_0^{+\infty}\hspace{-0.5cm}e^{-u}u^{i\beta}du+O\left(\lambda(s)+\frac{1}{s^2}\right)\right] \cr\hfill
  =  \frac{e^{i\frac{s^2}{2}\alpha}e^{-2i\beta\ln|s|}}{s^2}e^{i\beta\ln2}e^{-\pi\beta/2}2i\Gamma(1+i\beta)w(\alpha)
 +O\left(\frac{\lambda(s)}{s^2}+\frac{1}{s^4}\right),\cr}$$
\begin{eqnarray*}
\int_D^E e^{i\frac{s^2}{2}\tau}w(\tau)\frac{d\tau}{\sqrt{\tau-\alpha}} & = & \frac{e^{i\frac{s^2}{2}\alpha}}{|s|}\sqrt{2}e^{i\pi/4}\int_{\lambda(s)/2}^{s^2\delta/2}e^{-u}u^{-1/2}w\left(\alpha+2\frac{i}{s^2}u\right)du \\
 & = & \frac{e^{i\frac{s^2}{2}\alpha}}{|s|}\sqrt{2}e^{i\pi/4}\left[w(\alpha)\int_0^{+\infty}e^{-u}u^{-1/2}du+O\left(\sqrt{\lambda(s)}+\frac{1}{s^2}\right)\right] \\
 & = & \frac{e^{i\frac{s^2}{2}\alpha}}{|s|}\sqrt{2}e^{i\pi/4}w(\alpha)\Gamma\left(\frac{1}{2}\right)+O\left(\frac{\sqrt{\lambda(s)}}{|s|}+\frac{1}{|s|^3}\right).
\end{eqnarray*}
The result then follows from the choice $\lambda(s)=1/s^2$ (respectively $\lambda(s)=1/s^4$) for the logarithm (respectively the square root) and the fact that $\Gamma(1/2)=\sqrt{\pi}$.
\end{enumerate}
\end{dem}

Let us now give explicitly the contributions in the case of the $\Lambda^\cap$ contour proceeding with integration by parts and several applications of formula \aref{logarithm}.
\begin{eqnarray}
v^\cap_+(s)  & = & 2i\frac{z}{|z|}e^{-\pi|z|^2/8}\frac{\sinh\left(\pi|z|^2/4\right)}{|z|^2/4}\Gamma\left(1-i\frac{|z|^2}{4}\right)e^{i\left[\frac{s^2}{2}+\frac{|z|^2}{2}\ln|s|\right]}+O\left(\frac{1}{s^2}\right) \label{asymcap+} \\
v^\cap_-(s)  & = & 2\frac{\bar{z}}{|z|}e^{-\pi|z|^2/8}\frac{\sinh\left(\pi|z|^2/4\right)}{|z|^2/4}\Gamma\left(1+i\frac{|z|^2}{4}\right)e^{-i\left[\frac{s^2}{2}+\frac{|z|^2}{2}\ln|s|\right]}+O\left(\frac{1}{s^2}\right) \label{asymcap-} \\
v^\cap_0(s) 
 & = & O\left(\frac{1}{|s|}\right). \label{asymcap0}
\end{eqnarray}

The asymptotics for the $\Lambda^\cup$ contour can be mostly deduce from the following considerations (up to exponentially small errors):
\begin{itemize}
\item
the contribution of the piece $CD$ is the same as $B'A'$ whereas the contribution of $AB$ is the opposite of $A'B'$ (as we have looped once around the origin, the square root is changed into its opposite),
\item
the contribution of $DA$ is $-e^{\pi|z|^2/2}$ times the contribution of $BC$ (as we have looped once around $1$, $\ln(\tau-1)$ is translated by $2i\pi$ and the orientation of the path is reversed),
\item
in the case of $v_\pm(s)$, the contribution of $BC$ is $O(1/|s|)$ whereas in the case of $v_0(s)$ it is $O(1)$ and given by \aref{squareroot}.
\end{itemize}

\begin{eqnarray}
v^\cup_+(s) & = & -2i\frac{z}{|z|}e^{-\pi|z|^2/8}\frac{\sinh\left(\pi|z|^2/4\right)}{|z|^2/4}\Gamma\left(1-i\frac{|z|^2}{4}\right)e^{i\left[\frac{s^2}{2}+\frac{|z|^2}{2}\ln|s|\right]}+O\left(\frac{1}{|s|}\right) \label{asymcup+} \\
v^\cup_-(s) & = & 2\frac{\bar{z}}{|z|}e^{-\pi|z|^2/8}\frac{\sinh\left(\pi|z|^2/4\right)}{|z|^2/4}\Gamma\left(1+i\frac{|z|^2}{4}\right)e^{-i\left[\frac{s^2}{2}+\frac{|z|^2}{2}\ln|s|\right]}+O\left(\frac{1}{|s|}\right) \label{asymcup-} \\
v^\cup_0(s) 
 & = & {\rm sgn}(s)\frac{8}{|z|}e^{i\pi/4}\sqrt{\pi}\sinh\left(\frac{\pi|z|^2}{4}\right)+O\left(\frac{1}{|s|}\right). \label{asymcup0}
\end{eqnarray}

For later use, we introduce the functions
$$
a(z)=2\frac{z}{|z|}e^{-\pi|z|^2/8}\frac{\sinh\left(\pi|z|^2/4\right)}{|z|^2/4}\Gamma\left(1-i\frac{|z|^2}{4}\right), \quad
b(z)=\frac{8}{|z|}e^{i\pi/4}\sqrt{\pi}\sinh\left(\frac{\pi|z|^2}{4}\right),
$$
$$
\varphi(s,z)=\frac{s^2}{2}+\frac{|z|^2}{2}\ln|s|.
$$
Using that for a real number $u$ we have $\overline{\Gamma(1+iu)}=\Gamma(1-iu)$ and $|\Gamma(1+iu)|^2=\frac{\pi u}{\sinh(\pi u)}$, we get the useful relations
\begin{equation}
|b|^2=(2e^{\pi|z|^2/2}-2)|a|^2, \qquad \sqrt{2|a|^2+|b|^2}=\sqrt{2}e^{\pi|z|^2/4}|a|.
\end{equation}

%%%%%%%%%%%%%%%%%%%%%%%%%%%%%%%%%%%%%%%%%%%%%%%%%%%%%%%%%%%%%%%%%%%%%%%%%%%%%%%%%%
\subsection{The third solution}

As an application of Theorem~\ref{third} of the Appendix, let $V^\cup(s)$ and $V^\cap(s)$ denote the two independent solutions built in the preceding section through formula \aref{solcontint} and the corresponding $\Lambda^\cup$ and $\Lambda^\cap$ contours represented in Figure \ref{contour1}. We set
\begin{equation}
V^{\cup\cap}(s):=\frac{V^\cup(s)}{\Vert V^\cup(s)\Vert}\wedge\frac{V^\cap(s)}{\Vert V^\cap(s)\Vert}
\end{equation}
and it follows from the preceding result that $V^{\cup\cap}(s)$ is a solution of \aref{PJTsystem}. Using \aref{coordinates}, the asymptotics of our three (orthonormal!) solutions are
$$
\frac{V^\cup(s)}{\Vert V^\cup(s)\Vert}=\frac{\sqrt{2}e^{-\pi|z|^2/4}}{2|a|}\left(\begin{array}{c}
-iae^{i\varphi} \\
\overline{a}e^{-i\varphi} \\
{\rm sgn}(s)b
\end{array}\right)+O\left(\frac{1}{|s|}\right), \quad
\frac{V^\cap(s)}{\Vert V^\cap(s)\Vert}=\frac{\sqrt{2}}{2|a|}\left(\begin{array}{c}
iae^{i\varphi} \\
\overline{a}e^{-i\varphi} \\
0
\end{array}\right)+O\left(\frac{1}{|s|}\right),
$$
$$
V^{\cup\cap}(s)=\frac{e^{-\pi|z|^2/4}}{2|a|^2}\left(\begin{array}{c}
-{\rm sgn}(s)a\overline{b}e^{i\varphi} \\
-{\rm sgn}(s)i\overline{a}\overline{b}e^{-i\varphi} \\
2i|a|^2
\end{array}\right)+O\left(\frac{1}{|s|}\right).
$$

%%%%%%%%%%%%%%%%%%%%%%%%%%%%%%%%%%%%%%%%%%%%%%%%%%%%%%%%%%
\subsection{Scattering matrix}

To summarize, we have
$$
\lim_{s\to\pm\infty}\left(\begin{array}{ccc}
e^{i\left[\frac{s^2}{2}+\frac{|z|^2}{2}\ln|s|\right]} & 0 & 0 \\
0 & e^{-i\left[\frac{s^2}{2}+\frac{|z|^2}{2}\ln|s|\right]} & 0 \\
0 & 0 & 1
\end{array}\right)\left(\begin{array}{ccc} \frac{V^\cup(s)}{\Vert V^\cup(s)\Vert} & \frac{V^\cap(s)}{\Vert V^\cap(s)\Vert} & V^{\cup\cap}(s)\end{array}\right)=M_\pm
$$
where
$$
M_\pm:=\left(\begin{array}{ccc}
-i\frac{\sqrt{2}}{2}\frac{a}{|a|}e^{-\pi|z|^2/4} & i\frac{\sqrt{2}}{2}\frac{a}{|a|} & \mp i\frac{\sqrt{2}}{2}\frac{a}{|a|}\frac{\overline{b}}{|b|}\sqrt{1-e^{-\pi|z|^2/2}} \\
\frac{\sqrt{2}}{2}\frac{\overline{a}}{|a|}e^{-\pi|z|^2/4} & \frac{\sqrt{2}}{2}\frac{\overline{a}}{|a|} & \mp i\frac{\sqrt{2}}{2}\frac{\overline{a}}{|a|}\frac{\overline{b}}{|b|}\sqrt{1-e^{-\pi|z|^2/2}} \\
\pm\frac{b}{|b|}\sqrt{1-e^{-\pi|z|^2/2}} & 0 & ie^{-\pi|z|^2/4}
\end{array}\right)
$$
and the scattering matrix is given by $S=M_+M_-^{-1}$ which closes the proof of Theorem~\ref{th:scat}.

%%%%%%%%%%%%%%%%%%%%%%%%%%%%%%%%%%%%%%%%%%%%%%%%%%%%%%%%%%%%%%%%%%%%%%%%%%%%%%%%%%%%%%%%%%%%%%%%%%%

\section{Two-scale Wigner measure and Pseudo Jahn-Teller
Hamiltonians}\label{sec:measures}

In this section, we shall prove Theorem~\ref{theo:measures}. We proceed in two steps: we first prove the propagation of two-scale Wigner measure at infinity, then we calculate the transfer's rates at finite distance thanks to the scattering result of the preceding section. We shall work in the coordinates $(s,z,\sigma,\zeta)$ introduced in Section~\ref{sec:nf} and then translate the result in the original coordinates $(q,t,p,\tau)$. For this reason, we   consider the family $(v^\eps)$ of Theorem~\ref{theo:nf}. It satisfies microlocally in $\tilde\Omega$,
$$\op_\eps(Q)v^\eps=\eps f^\eps,$$
where for any $\chi\in{\cal C}_0^\infty (\tilde \Omega)$, $\op_\eps(\chi)f^\eps$ is uniformly bounded in $L^2(\R^3)$ and where $Q$ is given in~\aref{def:Q}.
For $\ell\in\{0,+1,-1\}$, we denote by 
$\tilde\Pi^\ell$  the eigenprojectors associated with the eigenvalues 
$$\tilde\lambda^\ell(s,z,\sigma,\zeta)=\sigma -\ell\,\sqrt{ s^2 +( \alpha(\zeta,\sigma))^2 z_1^2}$$ of $Q$.  The choice of the labelling is made so that $J^\pm_{in}\cup J^\mp_{out}\subset \{\tilde \lambda^\pm=0\}.$ On~$I$, one has $z_1=0$, therefore if $(e_1,e_2,e_3)$ is the canonical basis of $\C^3$,
\begin{equation}\label{proj}
\tilde \Pi^0=e_3\otimes e_3,\;\;\;\;
 \tilde\Pi^-=\left\{ \begin{array}{l} e_2\otimes e_2 \;{\rm if }\; s>0,
 \\ e_1\otimes e_1
\;{\rm if} \;s<0,
\end{array}\right.\;\;\;{\rm and}
\;\;\;\tilde\Pi^+=\left\{ \begin{array}{l} e_1\otimes e_1 \;{\rm if }\; s>0,\\ e_2\otimes e_2\;{\rm if} \;s<0.\end{array}\right.
\end{equation}
Let us now investigate
the implications of the change of coordinates on the two-scale Wigner measures.

$ $

\ni Let $\nu$ be the two-scale Wigner measure associated with the concentration at the scale $\sqrt\eps$ of $(\psi^\eps)$ on $I=\{q\wedge p=0\}$ with the choice of coordinates on $N(I)$ of Theorem~\ref{theo:measures}. We recall that these coordinates $(q,t,p,\tau,\eta)$ are such that for $\rho=(q,t,p,\tau)\in I$, the class on $N_\rho(I)$ of  the vector  $\delta\rho=(\delta q,\delta t ,\delta p,\delta \tau)$ is characterized by the function
$$\eta=\delta q\wedge p.$$
Let us denote by $N(\kappa)$ the change of coordinates induced on $N(I)$ by $\kappa$. 
By the geometric properties of two-scale Wigner measures and Egorov Theorem (see \cite{FG02}), the measure~$\tilde\nu=\nu\circ N(\kappa)$ is the two scale Wigner measure of $(v^\eps)$ for $I$ affected with the equation 
\begin{equation}\label{def:f}
g(s,z,\sigma,\zeta) z_1=0\;\;{\rm where}\;\;g\circ \kappa^{-1}(q,t,p,\tau)={|p|\over \lambda(p,\tau)\mu(p,\tau)}.
\end{equation}
Besides, $\tilde \nu$ decomposes as 
\begin{equation}\label{nudecomp}
\tilde\nu=\tilde \nu^0 \,\tilde\Pi^0 + \tilde \nu^+ \,\tilde \Pi^++\tilde \nu^-\,\tilde \Pi^-
\end{equation}
and we are concerned with the link  of the traces $\tilde\nu^\ell_{out}$ of $\tilde\nu^\ell$ on $S$ for $s>0$ with the traces $\tilde\nu^\ell_{in}$ of $\tilde\nu^\ell$ on $S$ for $s<0$.

$ $

\ni Finally, let us point out that  it is not clear that the system~\aref{systreduit} is well posed in~$s$ (indeed, the function $\alpha$ depends on the variable $\sigma$). We address this point in the first subsection, then we successively analyze the two scale Wigner measure above $|\eta|<\infty$ and $|\eta|=\infty$ in the two following subsections.

\subsection{An hyperbolic estimate}

\begin{prop}
\label{prop:esthyp}
Let $( v^\eps)$ be a bounded family of  
$L^2(\R_s\times\R_z^2,{\bf C}^3)$ 
satisfying \emph{\aref{systreduit}}, then 
if $\phi\in{\cal C}_0^\infty(\R)$ with $\phi=1$ close to $0$, there exists 
   $\delta_0>0$  and $\eps_0>0$ such that the family 
   $\bigl(\phi({z_1\over\delta_0}) v^\eps\bigr)_{\eps<\eps_0}$
  is bounded in  $L^\infty\bigl(\R_s,L^2(\R^2_z)\bigr)$.
\end{prop} 

\ni This proposition is proved in the same manner as Proposition~7 in \cite{FG02} and implies as in \cite{FG02} that one can evaluate the action of pseudodifferential operators on functions $(v^\eps)$ as follows.

\begin{lem}\label{lem:est}
There exists a neighborhood  $ \Omega_1\subset \tilde\Omega$ of  $0$ such that for all $\chi\in{\cal C}_0^\infty( \Omega_1)$,
the family $\underline v^\eps=\op_\eps(\chi)v^\eps$ 
satisfies:\\
 for $\eps$ small enough and for all compact ${\mathcal K}\subset \R^6$, there exists $N\in\N$ such that for all 
 $a\in{\cal C}_0^\infty({\mathcal K},\C^{3,3})$,
   $$\Bigg|\Bigl(\op_\eps(a) \underline v^\eps\;\mid\;\underline v^\eps\Bigr)\Biggr|
   \leq C\int_{-\infty}^{+\infty}
   \build{\rm sup}_{k+\mid\beta\mid\leq N}^{}
   \build{\rm sup}_{(z,\sigma,\zeta)\in\R^5}^{}
   \mid\partial_\sigma^k\partial_\zeta^\beta a(s,z,\sigma,\zeta)\mid \,
   ds.$$
\end{lem}

\ni We refer to  \cite{FG02} for the proofs. For the proof of the Proposition, one just need to replace the matrix~$L$ of~\cite{FG02} by the new definition \aref{eq:L} and for the proof of the Lemma, one has to exchange the role of $z$ and $\zeta$. 

\begin{rmq}\label{rem:est} Note that if the symbol of Lemma~\ref{lem:est} is of the form $a(s,z,\sigma,\zeta)=a_0(z,\sigma,\zeta)\phi\left({s\over r}\right)$ with~$r>0$, $\phi\in{\cal C}_0^\infty(\R)$ and $a_0\in{\cal C}^\infty_0(\R^5)$, then 
$$\left|\Bigl(\op_\eps(a) \underline v^\eps\;\mid\;\underline v^\eps\Bigr)\right|
   \leq C\int_{-\infty}^{+\infty} \left|\phi\left({s\over r}\right)\right| d s\td_r,0 0.$$
   Therefore, if $\tilde\mu$ is a Wigner measure of $(v^\eps)$, we have $\tilde\mu\left(\{s=0\}\right)=0$.  
\end{rmq}

\subsection{Analysis at finite distance} \label{subsec:finitedist}

We aim at proving Theorem~\ref{theo:measures} in $\{|\eta|<+\infty\}$, i.e. in  $\{|\eta|<R\}$ for any $R>0$. Let us fix $R>0$, we use the two-microlocal normal form of Section~\ref{sec:2micro}:
 we localize $\eta$ in the ball $B=\{|\eta|\leq R\}$ of $\R$ applying a cutoff function,
 $\phi({\eta\over R})$
for $R>0$ and  $\phi\in{\cal C}_{0}^\infty(\R)$.
 Then,  applying Proposition ~\ref{prop:nf2}, we set $V^\eps=\op_\eps^I({\rm Id}+\sqrt \eps \tilde C_1)v^\eps$, and we
 have
 \begin{equation}\label{system:codim3}
     \forall a\in{\cal
 C}_{0}^\infty(\R^6\times B),\;\;\op_{\eps}^I(a)\op_{\eps}( Q_0) V^\eps=\eps F^\eps(s,z)
 \end{equation}
 where $Q_0$ is defined in \aref{def:Q0} and $(F^\eps)$ is uniformly bounded
 in $L^2(\R^3_{s,z})$.

\noindent  Obviously, the
 two-scale Wigner measures of $(V^\eps)$ and $(v^\eps)$ are the same in $B$ and 
we   aim to describe the traces on $s=0^{-}$ in terms of the traces
on $s=0^+$ of the two--scale semiclassical measures $\tilde \nu$ of
$(V^\eps)$. The function
 $(F^\eps)$ does not contribute to this description.
Indeed,
denote by $S_{\eps}(s,s')$ the evolution operator associated with the
free system \aref{eq:2}
 and let $(U^\eps)$ be a solution of this system defined more precisely by,
\begin{equation}\label{eq:2}
\op_{\eps}(Q_0)U^\eps=0,\;\;U^\eps_{|s=0}=V^\eps_{|s=0},
\end{equation}
then we have
$$V^\eps(s)=U^\eps(s)+i\int_{0}^{s} S_\eps(0,t) F^\eps(t)dt.$$
Since the family $(F^\eps)$ is uniformly bounded in $L^{2}(\R^3_{s,z})$, we
deduce from H\"older inequality,
$$V^\eps(s)=U^\eps(s)+O(\sqrt{\mid s\mid})\;\;{\rm in}\;\;L^2(\R^{d}_{z}).$$
Therefore, the traces of the two-scale Wigner measure of $(U^\eps)$ and $(V^\eps)$ on $s=0^\pm$ are the same.
In the following, we focus on system \aref{eq:2}, which allows to calculate
the transfer coefficient in the variables $(s,z,\sigma,\zeta)$ and
then, in the variables $(q,t,p,\tau)$.

$ $

\noindent Let us 
  denote as before by $\tilde \nu$ the two-scale Wigner measure of the family $(v^\eps)$ (and thus of $(U^\eps)$ on~$I\times B$) for  the hypersurface $I=\{g(s,z,\sigma,\zeta)z_1=0\}$ where $g$ is defined in~\aref{def:f}.   Because of~\aref{proj} and~\aref{nudecomp}, one can   identify the measures $\tilde \nu^0$, $\tilde\nu^+$ and $\tilde\nu^-$ with the two-scale Wigner measures of the components of $U^\eps$. Set $U^\eps=(U_1^\eps, U_2^\eps, U_3^\eps)$, then $\tilde \nu^0$ is the two-scale
Wigner measure of $(U_3^\eps)$, $\tilde \nu^+$ of $(U_1^\eps)$ if $s>0$  and of $(U^\eps_2)$ if $s<0$ and finally, $\tilde \nu^-$  the measure of $(U^\eps_2)$ if $s>0$ and of $(U^\eps_1)$ if $s<0$.

$ $

\ni The crucial observation is that if one performs a symplectic  change of variables $(z,\zeta)\mapsto (\tilde z,\tilde\zeta) $ such that $\tilde z_1= \alpha(0,\zeta)z_1$ and if one sets 
$$\tilde U^\eps(s,\tilde z)=  K U^\eps(\sqrt \eps\, s,\sqrt\eps\, \tilde z)$$
where $K$ is a (scalar) Fourier integral operator associated with the symplectic change of coordinates,
then $\tilde U^\eps=(V_+,V_-,V_0)$ satisfies  the system \aref{PJTsystem} with $z=\tilde z_1$.
Then, the measure $\tilde \nu^+$ in $\{s>0\}$ is the two scale Wigner measure of $(V_+)$. 
By the scattering Theorem~\ref{th:scat} and the relation between  $\omega_+$ in one side and $\alpha_0$, $\alpha_+$ and $\alpha_-$ on the other side,  we obtain 
$$\nu^+_{out}=(1-{\rm e}^{-\pi \tilde\eta ^2/2})^2\,\nu^+_{in} + {\rm e}^{-\pi \tilde\eta^2}\;\nu^-_{in} + 2 (1-{\rm e}^{-\pi \tilde\eta^2/2}){\rm e}^{-\pi \tilde\eta^2/2}\, \nu^0_{in}$$
on the condition that $\nu^0_{in}$, $\nu^+_{in}$ and $\nu^-_{in}$ are two by two mutually singular and where the 
variable $\tilde \eta$ characterize the class of the vector $\delta\rho$ in $N_\rho(I)$ according to 
$$\tilde \eta=\delta \tilde z_1.$$
Therefore,
$$\tilde \eta= \alpha(0,\zeta) {\eta\over g(s,z,\sigma,\zeta)}.$$
In view of \aref{def:f}, $|\lambda|=|p|^{-1/2}$,  $\alpha(\sigma,\zeta)=\mu(\sigma,\zeta)^{-1}$ and of the fact that $\sigma=0$ on $I$, we obtain 
$$\tilde\eta=\eta\,|p|^{-{3\over 2}}.$$
Studying similarly $\nu^-_{out}$ and $\nu^0_{out}$, we obtain Theorem~\ref{theo:measures} in $|\eta|<+\infty$.

\subsection{Analysis at infinity}

It remains to prove Theorem~\ref{theo:measures} in $\{|\eta|=\infty\}$, i.e. that
$$\tilde\nu ^\ell_{out} {\bf 1}_{|\eta|=\infty} = \tilde\nu ^\ell_{in} {\bf 1}_{|\eta|=\infty}\;\;\forall\ell\in\{0,+1,-1\}.$$
Let us argue for the $+$ mode.  We follow the strategy of \cite{FG02} in Section 7.2 and detail the arguments because the fact that we have three modes implies a few complications.  Remark~\ref{rem:est} shows that $\tilde\nu^+(\{s=0\})=0$.
We
 consider $r>0$, $a_0=a_0(z,\sigma,\zeta)\in{\cal
C}_0^\infty (\R^5)$, $\phi\in{\cal
C}_0^\infty(\R)$, $\phi= 1$ in a neighborhood of  $0$,
$\delta<\delta_0$ (where $\delta_0$ is given by Proposition
\ref{prop:esthyp}), $R>0$.
We suppose $r$ small enough and ${\rm supp}(a_0)$ conveniently chosen so that $a(\cdot,\eta)$ is compactly supported in $\Omega_1$ for all $\eta\in\R$, where~$\tilde\Omega_1$ is the open set of Lemma~\ref{lem:est}.
 Then, we have
$$\left<a_0\;,\;\tilde\nu ^+_{out} {\bf 1}_{|\eta|=\infty} -\tilde\nu ^-_{in} {\bf 1}_{|\eta|=\infty}\right>= \lim_{r\rightarrow 0}
\left<a_0 \phi\left({s\over r}\right)\;,\;\partial_s \tilde \nu^+{\bf 1}_{|\eta|=+\infty}\right>.$$
We set
$$a(s,\sigma,z,\zeta,\eta)=a_0(z,\sigma,\zeta)
\phi\left({s\over r}\right)\phi\left({z_1\over\delta}\right)
\left(1-\phi\right)\left({\eta\over R}\right).$$
The function $a$ is a symbol of ${\cal A}$ and on the support of $a(s,\sigma,z,\zeta,z_1/\sqrt\eps)$, one has $|z_1|>R\sqrt\eps$ so that  
 the function
$a \tilde\Pi^+$
 is smooth, thus
 the operator  $\op_\eps^I(a \tilde\Pi^+)$ is well defined. Besides,
 $$\left<a_0 \phi\left({s\over r}\right)\;,\;\partial_s \tilde \nu^+{\bf 1}_{|\eta|=+\infty}\right>=\lim_{R\rightarrow+\infty} \,\lim_{\delta\rightarrow 0}\,\lim_{\eps\rightarrow 0}\left(\op_\eps^I(\partial_s a\tilde \Pi^+)v^\eps\;,\; v^\eps\right).$$
  Therefore, we set 
  $$L_\eps= \left(\op_\eps^I(\partial_sa\tilde \Pi^+)v^\eps\;,\; v^\eps\right)$$
  and we shall  study $L_\eps$ as first $\eps \rightarrow 0$,  then $\delta\rightarrow 0$, $R\rightarrow +\infty$ and finally $r\rightarrow 0$.
  
  $ $

  \noindent  We  remark that 
$$\op_\eps^I(\partial_sa\tilde \Pi^+) = {i\over\eps}\,\left[ {\eps\over i}\,\partial_s\;,\;\op_\eps^I(a\tilde\Pi^+)\right]-\op_\eps^I(a\partial_s\tilde\Pi^+).$$
  Therefore, using the equation,  we can  write 
  $L_\eps=L_\eps^1+L_\eps^2+L_\eps^3$ with 
  \begin{eqnarray*}
  L_\eps^1 & = & -{i}\left(\op_\eps^I(a\tilde\Pi^+)v^\eps,f^\eps\right) +i\left( \op_\eps^I(a\tilde\Pi^+) f^\eps,v^\eps\right)\\
  L_\eps^2 & = & {i\over \eps} \left(\left[ \op_\eps\left(
  V_{PJT}(s,\alpha(\sigma,\zeta)z_1)\right)
  \;,\;\op_\eps^I(a\tilde\Pi^+)\right]v^\eps,v^\eps\right)\\
  L_\eps^3 & = & \left( \op_\eps^I(a\partial_s\tilde\Pi^+) v^\eps,v^\eps\right). 
  \end{eqnarray*}
 We now prove successively that these three terms goes to $0$ in the prescribed limit.
 We will use symbolic calculus and for that, we introduce some notations. If 
 $q=q(s,z,\sigma , \zeta ,\eta )$, we set 
 $$\displaylines{
 q_\eps(s,z,\sigma,\zeta)= q\left(s,z,\sigma,\zeta,{z_1\over\sqrt\eps}\right)\;\;{\rm and}\;\;
 q_\eps^\sharp (s,z,\sigma,\zeta)= q\left(s,z,\eps\sigma,\eps\zeta,{z_1\over \sqrt\eps}\right)\cr}$$
 so that we have  
$$\op_\eps^{I}(q)= \;\op_\eps(q_\eps)=\;\op_1(q_\eps^\sharp).$$

\ni We observe that $\tilde \Pi^+ $ is homogeneous in the variables $(s,\alpha(\sigma,\zeta)z_1)$, thus for $q=a\tilde \Pi^+$, we have 
$$\mid \partial ^\beta _{\sigma ,z_2,\zeta }\partial ^\gamma
_{s,z_1}  q_\eps
(s,z,\sigma ,\zeta )\mid \leq C_{\beta ,\gamma }\left ({1\over R\sqrt
\eps }\right )^
{\mid \gamma \mid}\ $$
where we have used that 
\begin{equation}\label{def:theta}
\theta(s,z_1,\sigma,\zeta):=\sqrt{s^2 +\alpha(\sigma,\zeta)^2z_1^2}\geq R\sqrt\eps
\end{equation} on the support of $q_\eps$ and 
where we implicitely  assume 
$R\sqrt \eps\leq 1$.
In terms of  Weyl--H\"ormander metric (ses   sections 18.4, 18.5, 18.6  of \cite {Ho3}), 
one has  $q_\eps^\sharp \in S(1,g_\eps)$, where $g_\eps$ is the metric
$$g_\eps=dz_2^2+{ds^2+dz_1^2\over R^2\eps}+\eps^2(d\sigma ^2+d\zeta ^2).$$
Then
$${g_\eps\over g_\eps^\sigma }\leq \left ({\sqrt \eps\over R}\right )^2\ ,$$
so that the gain of the symbolic calculus is of order 
$\displaystyle {{\sqrt \eps\over R}}$.

$ $

 $\bullet$ {\bf Analysis  of $L_\eps^1$}:
 Since the function $f^\eps$ is locally bounded in $L^2$, we consider $\chi\in{\mathcal C}_0^\infty(\Omega_1)$ such that $\chi =1$ on the support of $a$ (as before the open set $\Omega_1$ is the one of Lemma~\ref{lem:est}). Then
 the symbolic calculus Theorem  18.5.4 of \cite{Ho3} gives
$$ \op_\eps^{I}(a\tilde\Pi ^+)= \, \op_\eps^{I}(a\tilde\Pi ^+\chi )\in  {\rm
op}_\eps^{I}(a\tilde\Pi ^+)
\op_\eps(\chi )+\op_1\left (S\left (\left ({\sqrt \eps\over
R}\right )^N,g_\eps\right )\right
)$$ for all  $N\in \N$. Then, the $L^2$ continuity Theorem
18.6.3 of \cite{Ho3} implies as $\eps\rightarrow 0$,
$$( \op_\eps ^{I}(a\tilde\Pi ^+)f^\eps\;\mid\;v ^\eps)=
( \op_\eps^{I}(a\tilde\Pi ^+)\op_\eps(\chi
)f^\eps\;\mid\;v ^\eps)+o(1) .$$
Lemma~\ref{lem:est} and the boundedness in $L^2$ of $\op_\eps(\chi
)f^\eps$ gives the existence of $C>0$ such that 
$$\left|( \op_\eps^{I}(a\tilde\Pi ^+)\op_\eps(\chi
)f^\eps\;\mid\;v ^\eps)\right|\leq  C\,\int_{-\infty}^{+\infty} \left|\phi\left({s\over r}\right)\right| \d s\td_{r},{0} 0.$$
One argues similarly with the other term.

$ $

  $\bullet$ {\bf Analysis of $L_\eps^2$}:  We first observe that  $\left(V_{PJT}(s,\alpha(\sigma,\zeta)z_1)\right)_\eps^\sharp \in S( \theta_\eps^\sharp, g_\eps)$  where $\theta$ is defined in \aref{def:theta} and that 
  $q_\eps^\sharp \in S((\theta_\eps^\sharp)^{-N},g_\eps)$  for all $N\in \N$ because $a$ is supported in a fixed compact set. Therefore, using that 
  $V_{PJT}(s,\alpha(\sigma,\zeta)z_1)$ and $\tilde\Pi^+$ commute,
  we obtain by symbolic calculus
  $$\displaylines{\qquad
  {1\over \eps}
  \left[ \op_1\left(\left(V_{PJT}(s,\alpha(\sigma,\zeta)z_1)\right)_\eps^\sharp\right)
  \;,\;\op_1(q_\eps^\sharp)\right]
  \in\;\op_1 \Biggl({1\over 2i\eps} 
  \left\{\left( V_{PJT}(s,\alpha(\sigma,\zeta)z_1)\right)_\eps^\sharp
  \;,\;(a\tilde\Pi^+)_\eps^\sharp\right\} 
  \hfill\cr\hfill
  -{1\over 2i\eps} \left\{ (a\tilde\Pi^+)_\eps^\sharp\;,\;
 \left(V_{PJT}(s,\alpha(\sigma,\zeta)z_1)\right)_\eps^\sharp\right\}\Biggr) +{1\over \eps}\op_1\left(S\left (\left ({\sqrt \eps\over R}\right )^2,
g_\eps\right )\right ).\qquad\cr}$$
It is at this very place that the proof differs  from the one of \cite{FG02} because of the presence of the three modes which induces the presence of additional terms. 
We observe that 
$$\displaylines{\qquad{1\over \eps} \left\{ \left( V_{PJT}(s,\alpha(\sigma,\zeta)z_1)\right)_\eps^\sharp
  \;,\;(a\tilde\Pi^+)_\eps^\sharp\right\} 
 =  \left(\left\{V_{PJT}(s,\alpha(\sigma,\zeta)z_1)\;,\;(a\tilde \Pi^+)\right\}_{s,z,\sigma,\zeta}\right)^\sharp_\eps\hfill\cr\hfill +{1\over R\sqrt\eps}\left(\partial_\eta aV_{PJT}(0,\partial_{\zeta_1}\alpha(\sigma,,\zeta)z_1)\tilde\Pi^+\right)_\eps^\sharp\qquad\cr\hfill
 =   \left(\left\{V_{PJT}(s,\alpha(\sigma,\zeta)z_1)\;,\;(a\tilde \Pi^+)\right\}_{s,z,\sigma,\zeta}\right)^\sharp_\eps +{1\over R}\left(\partial_\eta a\,V_{PJT}(0,\partial_{\zeta_1}\alpha(\sigma,\zeta)\eta)\,\tilde\Pi^+\right)_\eps^\sharp.\qquad
\cr}$$
Similarly,
$$\displaylines{\qquad
{1\over \eps} \left\{ (a\tilde\Pi^+)_\eps^\sharp\;,\;
 \left(V_{PJT}(s,\alpha(\sigma,\zeta)z_1)\right)_\eps^\sharp\right\}
 \hfill\cr\hfill=\left(\left\{(a\tilde \Pi^+)\;,\;V_{PJT}(s,\alpha(\sigma,\zeta)z_1)\right\}_{s,z,\sigma,\zeta}\right)^\sharp_\eps +{1\over R}\left(\partial_\eta a\,\tilde\Pi^+\,V_{PJT}(0,\partial_{\zeta_1}\alpha(\sigma,\zeta)\eta)\right)_\eps^\sharp.\qquad
\cr}$$
 Therefore $L^2_\eps=\left(\op_\eps^I (B) v^\eps\;,\;v^\eps\right)+O(R^{-2})$ with
$$\displaylines{\qquad B=\left\{V_{PJT}(s,\alpha(\sigma,\zeta)z_1)
  \;,\;a\tilde\Pi^+\right\} _{s,z,\sigma,\zeta}
 - \left\{ a\tilde\Pi^+\;,\;
 V_{PJT}(s,\alpha(\sigma,\zeta)z_1)\right\}_{s,z,\sigma,\zeta} \hfill\cr\hfill
 + {1\over R}\partial_\eta a\left(\tilde\Pi^+\,V_{PJT}(0,\partial_{\zeta_1}\alpha(\sigma,\zeta)\eta)+V_{PJT}(0,\partial_{\zeta_1}\alpha(\sigma,\zeta)\eta)\,\tilde\Pi^+\right).\qquad\cr}$$
 We write  $B=B_1+B_2+B_3$ with 
$$B_1=\{V_{PJT},a\}\tilde\Pi^+ -\tilde\Pi^+\{a ,V_{PJT}\}\;\;{\rm and}\;\;
B_2=a\left(\{V_{PJT},\tilde\Pi^+\}-\{\tilde\Pi^+,V_{PJT}\}\right).$$
The matrix $B_1+B_3$ is such that by Lemma~\ref{lem:est}, we have 
$$\left|\left( \op_\eps^I(B_1+B_3) v^\eps,v^\eps\right) \right|\leq 
C\int _{-\infty
}^{+\infty }\left({\delta\over r }\left \vert \phi '\left (s\over r
\right )\right \vert +\phi
\left ({s\over r }
\right )\right )\, ds.$$
The term in $B_2$ is more complicated and we will cut it into two parts. Indeed, using
$$V_{PJT}(s,\alpha(\sigma,\zeta)z_1)=\theta\,(\tilde\Pi^+-\tilde\Pi^-)$$
we write 
\begin{eqnarray*}
B_2 & = & a\theta\left(\{\tilde \Pi^+-\tilde\Pi^-,\tilde\Pi^+\}-\{\tilde\Pi^+,\tilde\Pi^+-\tilde\Pi^-\}\right) + 
a\left((\tilde\Pi^+-\tilde\Pi^-) \{\theta,\tilde\Pi^+\}-\{\tilde\Pi^+,\theta\}(\tilde\Pi^+-\tilde\Pi^-)\right)\\
& = & a\theta\left(\{\tilde\Pi^+,\tilde\Pi^-\}-\{\tilde\Pi^-,\tilde\Pi^+\}\right) + a\left((\tilde\Pi^+-\tilde\Pi^-) \{\theta,\tilde\Pi^+\}+\{\theta,\tilde\Pi^+\}(\tilde\Pi^+-\tilde\Pi^-)\right).
\end{eqnarray*} 
We  shall prove that $B_2$ is  off-diagonal  where we say that the matrix $D$ is off-diagonal if $\Pi^\ell D\Pi^\ell =0$ for any $\ell\in\{0,+1,-1\}$.  We will study those terms together with $L^3_\eps$ since $\partial_s \tilde\Pi^+$  is also off-diagonal.

$ $

\ni Set $B_2=B_{2,1}+B_{2,2}$ with 
$$B_{2,2}=a\theta\left(\{\tilde\Pi^+,\tilde\Pi^-\}-\{\tilde\Pi^-,\tilde\Pi^+\}\right).$$
We remark that the fact  that $\tilde\Pi^\ell\tilde\Pi^{\ell'}=\delta_{\ell,\ell'}\tilde\Pi^\ell$ for all $\ell,\ell'\in\{0,+1,-1\}$ yields 
\begin{equation}\label{33'}\forall\ell\in\{0,+1,-1\},\;\;
\forall h\in{\cal C}^\infty(\R^6),\;\;\Pi^\ell \{h,\Pi^\pm\}\Pi^\ell=0.
\end{equation}
Therefore, applying \aref{33'} to $h=\theta$,  the matrix $B_{2,1}$ is of the form $B_2=aD$ with $D$ off-diagonal. 

$ $

\ni Let us now consider $B_{2,2}$. 
Using the general fact $A\{B,C\}-\{A,B\}C=\{AB,C\}-\{A,BC\}$, we observe that for $\ell\not=k$, $\ell,k\in\{0,+1,-1\}$
\begin{eqnarray}
\{\tilde \Pi^\ell,\tilde \Pi^k\}&=&\{(\tilde \Pi^\ell)^2,\tilde \Pi^k\}=\tilde \Pi^\ell\{\tilde \Pi^\ell,\tilde \Pi^k\}-\{\tilde \Pi^\ell,\tilde \Pi ^\ell\}\tilde \Pi^k,\label{crochet-kj}\\
-\{\tilde \Pi^k, \tilde \Pi^\ell\}&=&-\{\tilde \Pi^k,(\tilde \Pi^\ell)^2\}=\tilde \Pi^k\{\tilde \Pi^\ell,\tilde \Pi^\ell\}-\{\tilde \Pi^k,\tilde \Pi^\ell\}\tilde \Pi^\ell,\nonumber
\end{eqnarray}
yielding $\tilde\Pi^0\{\tilde\Pi^+,\tilde\Pi^-\}\tilde\Pi ^0= -\tilde\Pi^0\{\tilde\Pi^-,\tilde\Pi^+\}\tilde\Pi^0=0$ and 
\begin{eqnarray*}
\tilde\Pi^-\{\tilde\Pi^+,\tilde\Pi^-\}\tilde\Pi ^-= -\tilde\Pi^-\{\tilde\Pi^+,\tilde\Pi^+\}\tilde\Pi^- &\; {\rm and}\;&
\tilde\Pi^+\{\tilde\Pi^+,\tilde\Pi^-\}\tilde\Pi^ += -\tilde\Pi^+\{\tilde\Pi^-,\tilde\Pi^-\}\tilde\Pi^+,\\
\tilde\Pi^-\{\tilde\Pi^-,\tilde\Pi^+\}\tilde\Pi ^-= -\tilde\Pi^-\{\tilde\Pi^+,\tilde\Pi^+\}\tilde\Pi^-&\; {\rm and}\;&
\tilde\Pi^+\{\tilde\Pi^-,\tilde\Pi^+\}\tilde\Pi ^+= -\tilde\Pi^+\{\tilde\Pi^-,\tilde\Pi^-\}\tilde\Pi^+,
\end{eqnarray*}
We obtain $\tilde\Pi^\ell B_{2,2}\tilde\Pi^\ell=0$ for all $\ell\in\{0,+1,-1\}$ and $B_{2,2}$ also is off-diagonal.

$ $

$\bullet$ {\bf Analysis of $L_\eps^3$}: Equation \aref{33'} for $h=\sigma$ implies that $\partial_s\tilde\Pi^+$ is off-diagonal. 
We now consider  terms of the form 
$$L^3_\eps=( \op_\eps^{I}(aD) v^\eps\;\mid\;v^\eps)\ $$
for some matrix $D$ homogeneous of degree $-1$ in the variables $(s,\alpha(\sigma,\zeta)z_1)$ and such that $\tilde\Pi^\ell D\tilde\Pi^\ell=0$ for all $\ell\in\{0,+1,-1\}$.
Without loss of generality we can suppose that $D=\tilde \Pi^\ell D\tilde \Pi^{\ell'}$ for some $\ell,\ell'\in\{0,+1,-1\}$ with $\ell\not=\ell'$. Then there exists a real number $c_{\ell,\ell'}$ ($c_{\ell,\ell'}\in\{\pm1,\pm2\}$) such that 
$$\left[ D,V_{PJT}\right]=c_{\ell,\ell'} \theta D.$$
We use this relation to write 
$$aD=\left[{aD\over c_{\ell,\ell'}\theta}\;,\;V_{PJT}\right]= \left[{aD\over c_{\ell,\ell'}\theta}\;,\;-\sigma{\rm Id}+V_{PJT}\right],$$
which allows to reuse the equation. Again by 
Weyl--H\"ormander calculus with the metric $g_\eps$, we obtain 
 $(V_{PJT})_\eps^\sharp \in S(\theta _\eps^\sharp, g_\eps)$ and
$$\left ({a D\over c_{\ell,\ell'} \theta}\right )^\sharp _\eps\in
S\left ({1\over
\sqrt \eps R \theta _\eps^\sharp } ,g_\eps\right )$$
(again because of 
 $\theta _\eps^\sharp \geq C\sqrt \eps R$ on the support of 
$a_\eps^\sharp $). Therefore, we have
$$\displaylines{\qquad\qquad
   \op_\eps^{I}(a D)\in\; \op_\eps(-\sigma + V_{PJT})\,\op_\eps^I\left
({aD \over c_{\ell,\ell'}\theta}\right ) - \op_\eps ^I\left ({aD
\over c_{\ell,\ell'} \theta}\right )
\op_\eps(-\sigma+V_{PJT})\hfill\cr\hfill
+{\eps\over i}\op_\eps^I\left(\partial_s\left({aD\over c_{\ell,\ell'}\theta}\right)\right)
+{\rm op }_1\left (S\left ({1\over R^2},g_\eps\right )\right )
\qquad\qquad\cr} $$
whence, using the equation for 
 $v^\eps$,
$$\displaylines{
\qquad L^3_\eps
=O\left ({1\over R^2}\right )
+\eps\left( \op_\eps^{I}
\left(  {aD \over c_{\ell,\ell'}\theta}\right )v^\eps\;\mid\;f^\eps\right)
-\eps\left ( \op_\eps^I
\left({aD \over c_{\ell,\ell'}\theta}\right )f^\eps\;\mid\;v^\eps\right)\hfill\cr\hfill
+{\eps\over i}\left(\op_\eps^{I}\left(\partial _s\left({aD \over c_{\ell,\ell'}\theta}\right )\right)
v^\eps\;\mid\;v^\eps\right)\ .\qquad\cr}$$
Arguing as before, one gets as  $\eps$ goes to $0$ then  $R$ to $+\infty $,
$$L_\eps^3=o(1)+O\left ({1\over R^2}\right )+{\eps\over i}\left(\op_\eps^{I}\left(\partial _s\left
({aD \over c_{\ell,\ell'}\theta}\right )\right)
v^\eps\;\mid\;v^\eps\right).$$
Observing that  $$\left|\partial ^\beta _{\sigma ,\zeta}\,\partial _s\left({aD \over c_{\ell,\ell'}\theta}\right)_\eps\right|
\leq {C_\beta \over (s^2+\eps R^2)^{3/2}}$$
we obtain by Lemma~\ref{lem:est}
$$\left| L^3_\eps\right| \leq o(1)+O\left ({1\over R^2}\right )+C\eps\int
_{-\infty}^{+\infty }
{ds\over (s^2+\eps R^2)^{3/2}}
\leq o(1)+O\left ({1\over R^2}\right )\ ,$$
which finishes the proof.
  
 %%%%%%%%%%%%%%%%%%%%%%%%%%%%%%%%%%%%%%%%%%%%%
 
 \appendix
 \section{ Appendix: Complex wedge product and solutions of ODEs}

We consider $(\Hh,\langle\cdot|\cdot\rangle)$  an Hilbert space of dimension $3$ over $\C$ and  $A(s)$  a continuous family of selfadjoint endomorphisms on $\Hh$. We consider the Schr\"odinger system
\begin{equation} \label{Schrod}
-i\partial_s V=A(s)V.
\end{equation}
Let $[\cdot,\cdot,\cdot]$ denote a non trivial alternating trilinear form on $\Hh$ and let $u_1$ and $u_2$ be two vectors of $\Hh$. We define the wedge product of  $u_1$ and $u_2$ with respect to $[\cdot,\cdot,\cdot]$ to be the unique vector $u_1\wedge u_2$ (through Riesz representation theorem) such that
\begin{equation}
\forall w\in\Hh, \qquad [u_1,u_2,w]=\langle u_1\wedge u_2|w\rangle.
\end{equation}
Then, we have the two following points:
\begin{enumerate}
\item
If $A$ is an endomorphism of $\Hh$, we have 
\begin{equation}\label{rmq}[Au_1,u_2,u_3]+[u_1,Au_2,u_3]+[u_1,u_2,Au_3]=(\Tr A)[u_1,u_2,u_3].
\end{equation}
\item
If $u$ (respectively $v$) has coordinates $(u_1,u_2,u_3)$ (respectively $(v_1,v_2,v_3)$) in some orthonormal basis $(e_1,e_2,e_3)$ of $\C^3$ then the coordinates of $u\wedge v$ are
\begin{equation}\label{coordinates}
\overline{[e_1,e_2,e_3](u_2v_3-u_3v_2,u_3v_1-u_1v_3,u_1v_2-u_2v_1)}.
\end{equation}
\end{enumerate}

\begin{theo}\label{third}
Let $V_1(s)$ and $V_2(s)$ be two solutions of the system \aref{Schrod}. Then
\begin{equation}
e^{i\int^s_{s_0}\Tr A(s')ds'}V_1(s)\wedge V_2(s)
\end{equation}
is also a solution of \aref{Schrod}. In particular, if $V_1$ and $V_2$ are linearly independent, we have a basis of solutions.
\end{theo}

\ni Note that 
this is in fact a corollary of Liouville theorem on Wronskian combined with the fact that the equation is norm preserving.

$ $

\begin{dem}
If $V_1$ and $V_2$ are linearly dependent, the wedge product is zero and is certainly a solution. Thus, we can assume that $V_1$ and $V_2$ are linearly independent so that $(V_1(s),V_2(s),V_1(s)\wedge V_2(s))$ is an instantaneous basis of $\Hh$. Removing the phase, it is enough to show that $W(s)=V_1(s)\wedge V_2(s)$ is a solution of
$$
-i\partial_s W=[A(s)-\Tr A(s){\rm Id}]W
$$
which we will show to hold in the preceding instantaneous basis. We first have
$$
\langle V_j(s)|(-i\partial_sW)(s)\rangle  =  \langle-i\partial_sV_j|W\rangle=  \langle AV_j|W\rangle = \langle V_j|AW\rangle = \langle V_j(s)|[A(s)-\Tr A(s){\rm Id}]W(s)\rangle
$$
where we have used that $V_j$ and $W$ are orthogonal, the equation and the self-adjointness of $A(s)$.
Similarly, using~\aref{rmq} and~\aref{coordinates}, we obtain
\begin{eqnarray*}
\langle W(s)|(-i\partial_sW)(s)\rangle & = &- \langle W|(-i\partial_sV_1)\wedge V_2\rangle-\langle W|V_1\wedge(-i\partial_sV_2)\rangle \\
 & = &- \langle W|AV_1\wedge V_2\rangle-\langle W|V_1\wedge AV_2\rangle  \\
 & = & -\overline{[AV_1,V_2,W]}-\overline{[V_1,AV_2,W]} \\
 & = & \overline{[V_1,V_2,AW]-(\Tr A)[V_1,V_2,W]} \\
 & = & \langle[A-\Tr A{\rm Id}]W|V_1\wedge V_2\rangle \\
 & = & \langle W(s)|[A(s)-\Tr A(s){\rm Id}]W(s)\rangle.
 \end{eqnarray*}
 \end{dem}

%%%%%%%%%%%%%%%%%%%%%%%%%%%%%%%%%%%%%%%%%%%%%%%%%%%%%%%%%%%%%%%%%%%%%%%%%%%%%%%%%%%%%%%%%
\vskip 1cm 
\noindent{\bf Acknowledgements} : The authors thank A. Joye, P.
G\'erard and C. Lasser for friendly and fruitful discussions.

%%%%%%%%%%%%%%%%%%%%%%%%%%%%%%%%%%%%%%%%%%%%%%%%%%

\ni Clotilde \textsc{Fermanian Kammerer}, Universit\'e Paris Est, 
UFR des Sciences et Technologie,
61, avenue du G\'en\'eral de Gaulle,
94010 Cr\'eteil Cedex, France.\\
{\tt clotilde.fermanian@univ-paris12.fr}

\medskip

\ni Vidian \textsc{Rousse}, Universit\'e Paris Est, 
UFR des Sciences et Technologie,
61, avenue du G\'en\'eral de Gaulle,
94010 Cr\'eteil Cedex, France.\\
{\tt vidian.rousse@univ-paris12.fr}

\end{document}